\documentclass[11pt, a4paper]{amsart}
 \pdfoutput=1
\usepackage{fullpage}
\usepackage{amscd}
\usepackage{amsmath}
\usepackage{amsthm}
\usepackage{amssymb}
\usepackage{mathdots}
\usepackage[colorlinks,linkcolor=blue,filecolor=blue,citecolor=blue,urlcolor=blue,bookmarks=false,hypertexnames]{hyperref}
\usepackage{tikz}
\usepackage{caption}
\usepackage{subcaption}
\DeclareCaptionType[within=section]{mat}[Matrix]
\usepackage{paralist}

\def\unitcube (#1,#2,#3);{
\draw[](#1,#2,#3)--(#1+1,#2,#3);
}

\usetikzlibrary{decorations.pathreplacing}
\tikzset{
   mybrace/.style={decorate,decoration={brace,aspect=#1}}
}
\newcommand\naturals{\mathbb{N}}

\newcommand\reals{\mathbb{R}}
\def\MC{{\mathcal C}}

\def\MI{{\mathcal I}}
\def\MP{{\mathcal P}}
\def\MQ{{\mathcal Q}}
\def\MR{{\mathcal R}}

\def\opn#1#2{\def#1{\operatorname{#2}}} 
\opn\det{det}
\opn\ini{in}
\opn\height{height}
\opn\rank{rank}
\opn\supp{supp}
\opn\gcd{gcd}
\opn\Min{Min}
\opn\sp{sp}

\makeatletter
\@addtoreset{equation}{section}
\def\theequation{\thesection.\@arabic \c@equation}
\makeatother

\theoremstyle{plain}

\newtheorem{theorem}[equation]{Theorem}
\newtheorem{corollary}[equation]{Corollary}

\newtheorem{Lemma}[equation]{Lemma}
\newtheorem{proposition}[equation]{Proposition}

\theoremstyle{definition}

\newtheorem{definition}[equation]{Definition}

\newtheorem{discussion}[equation]{Discussion}
\newenvironment{discussionbox}[1][]{%
\begin{discussion}[#1]\pushQED{\qed}}{\popQED \end{discussion}}
\newtheorem{notation}[equation]{Notation}

\newtheorem{Remark}[equation]{Remark}
\newenvironment{remarkbox}[1][]{%
\begin{Remark}[#1]\pushQED{\qed}}{\popQED \end{Remark}}
\newtheorem{example}[equation]{Example}
\newenvironment{examplebox}[1][]{%
\begin{example}[#1]\pushQED{\qed}}{\popQED \end{example}}
\newcommand{\define}[1]{\emph{#1}}

\newtheorem{setup}[equation]{Setup}
\newtheorem{lemma}[equation]{Lemma}
\usepackage{enumitem}
\def\sgn{{\operatorname{sgn}}}

\title {Polyominoes and Knutson ideals}

\author{Mitra Koley}
\address{School of Mathematics, Indian Institute of Science Education and Research Thiruvananthapuram, Vithura, 695551, India}
\email{mitra@iisertvm.ac.in}

\author{Nirmal Kotal}
\address{The Institute of Mathematical Sciences, A CI of Homi Bhabha National Institute, Taramani, Chennai 600113, India}
\email{nirmalk@imsc.res.in}

\author{Dharm Veer}
\address{Department of Mathematics and Statistics, Dalhousie University, Halifax, B3H 4R2, Canada }
\email{d.veer@dal.ca}

\thanks{NK was supported by the Institute of Mathematical Sciences postdoctoral fellowship. DV was supported by AARMS postdoctoral fellowship.}



\begin{document}

\begin{abstract}
In this article, we study two fundamental questions on polyomino ideals which are radicality and primality. In order to study the question of radicality, we initiate the study of Knutson ideals among polyominoes. Knutson ideals were introduced by Conca and Varbaro after the work of Knutson on compatibly split ideals. 
Knutson ideals are known to have nice properties, for example, they are well behaved with Gr\"{o}bner bases, and it has square-free initial ideals; hence they are radical. 

We show that polyomino ideals associated with closed path, weakly closed path, simple thin, and ladder polyominoes are Knutson. We also show that polyomino ideals associated with a class of thin polyominoes are Knutson; hence they are radical. In fact, we show that these polyomino ideals are prime and the reduced Gr\"{o}bner basis is computed. Furthermore, we prove that under a certain condition, if a parallelogram polyomino is extracted from another parallelogram polyomino, the resulting collection of cells is Knutson. We also compute their Gr\"{o}bner basis.
\end{abstract}

\maketitle

\section{Introduction}
Polyominoes are plane figures obtained by joining unit squares (cells) edge to edge. They were formally defined by Golomb \cite{Golomb_polyomino, Golomb}. Polyominoes
appear in many recreational mathematics and combinatorics for instance, tiling  problem of a certain region or the plane with polyominoes.

In connection with commutative algebra, Qureshi~\cite{QUR12} assigns to every polyomino $\mathcal{P}$ the ideal of inner $2$-minors, called a polyomino ideal, denoted by $I_{\mathcal{P}}$. The class of polyomino ideals generalizes many known classes of ideals such as the ideal generated by  $2$-minors of a  generic matrix, the ideal generated by $2$-minors of one-sided or two-sided ladders, join-meet ideal of a planar distributive lattice. Determinantal ideals and ladder determinanatal ideals are well studied in past \cite{Conca-1, Conca-2, Conca-3, Gorla}. 

Let $X$ be a generic matrix of size $m\times n$. Let $I$ be an ideal generated by an arbitrary set of $2$-minors of $X$, then the natural question arises when $I$ is radical or prime or more generally what is its primary decomposition.
The motivation for studying ideals generated by arbitrary sets of $2$-minors of a generic matrix comes from algebraic statistics \cite{Algebraic-Statistics, Sturmfels-statistics}: more specifically, a primary decomposition of such ideals helps to measure the connectedness of the set of $m \times n$ contingency tables with the same row and column sums, see \cite{Eisenbud-Sturmfels}. In \cite{binomial-edge-ideal} the author showed that  if $X$ is a $2 \times n$ matrix and $I$ is an ideal generated by an arbitrary set of $2$-minors of $X$, then $I$ is always radical, in fact it is shown that $I$ has a squarefree initial ideal with respect to graded lexicographic ordering. Moreover the authors also gave a primary decomposition of $I$. The problem becomes much
more challenging if $m, n \geq 3$, as easy examples show that $I$ need not to be radical in general. The class of ideals of adjacent $2$-minors are studied in \cite{Hosten-Sullivant, Herzog-Hibi, Ohsugi-Hibi}.

Given the preceding, two interesting questions on polyomino ideals are radicality and primality of polyomino ideals. 
While there are no known examples of polyomino ideals that are not radical, there are polyominoes for which the associated ideal is not prime.
We say that a polyomino $\MP$ is {\em prime} (resp.~{\em non-prime}) if the polyomino ideal $I_\MP$ is prime (resp.~non-prime).
The question of primality of polyomino ideals were first studied by Qureshi~\cite{QUR12}, where she proved that convex polyominoes are prime.
Herzog and Madani~\cite{HM14} extended Qureshi's result and established that simple polyominoes are prime.
Qureshi et al.~\cite{QSS17} showed polyomino ideal of a simple polyomino
coincides with the toric ideal of a weakly chordal bipartite graph.
Additionally, Hibi-Qureshi~\cite{HQ15nonsimple} and Shikama~\cite{SHI18nonsimple} showed that non-simple polyominoes formed by removing a convex polyomino from a rectangle are prime.
Characterizing prime polyominoes is a complex problem, and to tackle it, Mascia et al.~\cite{MRR20primalityOfPolyominoes} introduced the notion of zig-zag walks for polyominoes. 
They showed~\cite[Corollary\ 3.6]{MRR20primalityOfPolyominoes} that the presence of zig-zag walks indicates that a polyomino is not prime, and they conjectured~\cite[Conjecture\ 4.6]{MRR20primalityOfPolyominoes} that this characteristic characterize non-prime polyominoes.
This conjecture has been confirmed for specific classes of polyominoes, 
including grid polyominoes~\cite{MRR20primalityOfPolyominoes}, closed path polyominoes~\cite{CN23closedpathprimality}, and weakly closed path polyominoes~\cite{CNU22weaklyclosedprimality}. 
Furthermore, Mascia et al.~\cite{MRR22:primality_of_polyomino} used Gr\"{o}bner bases to identify certain  prime polyominoes.
The only known classes of non-prime radical polyomino ideals are closed path polyominoes and weakly closed path polyominoes with zig-zag walks.
In \cite{CNV2024polyocollection}, Cisto, Navarra and the last author of this paper studied the primary decomposition of non-prime closed path polyomino ideals and established that these polyomino ideals are unmixed.

In recent years, various algebraic properties of algebras associated with polyominoes have been explored by several researchers. 
For instance, the Hilbert series of these algebras has been examined in \cite{RR21,KV21,QRR22,KV22,JN24rookpolynomial,CNU2022hilbert}, Cohen-Macaulayness has been investigated in \cite{QUR12, HM14, CNU22grobnerclosedpath}, and the Gorenstein condition has been studied in \cite{QUR12, EHQR21,RR21, QRR22, CNU2022hilbert}. 
Interested readers may refer to \cite{CNV2024polyocollection} for a generalization of polyominoes.

 In this article, in order to study radicality of polyomino ideals, we initiate to study Knutson ideals among polyominoes. Knutson ideals were introduced by Conca and Varbaro in \cite{Conca-Varbaro}, in prime characteristic which corresponds to compatibly split ideals. In arbitrary characteristic, Knutson ideals have squarefree initial ideals, hence the are radical. Recently it has already been shown that some interesting classes of ideals of minors are Knutson ideals; for instance, determinantal ideals of Hankel matrices \cite{Lisa1}, generic matrices \cite{Lisa2}, ideals generated by arbitrary set of $2$-minors in a generic $2\times n$ matrix \cite{Lisa-binomial}, ladder determinantal ideals \cite{DS-NB-S-V}. We say that a polyomino $\MP$ is {\em Knutson} if the polyomino ideal $I_\MP$ is a Knutson ideal.

 We show that simple thin polyominoes, closed path polyominoes, weakly closed path polyominoes, ladder polyominoes, various thin collection of cells are Knutson, hence they are radical. We also compute Gr\"obner bases for some of these classes of collection of cells. As an intermediate step of showing Knutson for simple thin polyominoes, closed path polyominoes, and weakly closed path polyominoes, we show that they are of K\"{o}nig type. Herzog et al.~\cite{HHM22konigtype} introduced the concept of ideals of K\"{o}nig type to examine the Cohen-Macaulay property of their initial ideals, independent of the characteristic of the base field. 
 In this article, we show that the notion of K\"{o}nig type can be employed to determine Knutson ideals.
 Herzog and Hibi~\cite{HH23konig} showed that simple thin polyominoes are of K\"{o}nig type. However, we give an example that indicates their argument was not exhaustive; a more general version of their pivotal lemma \cite[Lemma\ 5.3]{HH23konig} is necessary (see Lemma~\ref{lem:herzoghibi}). Utilizing this refined lemma, we complete the proof of  the theorem which states that simple thin polyominoes are indeed of K\"{o}nig type. Dinu and Navarra~\cite{DN22konigtype} showed that closed path polyominoes are of K\"{o}nig type. We show that weakly closed path polyominoes are of K\"{o}nig type. Recently, this result was also outlined by Navarra in \cite{navarra2024weaklyclosedkonig}.
 
 We then proceed to show that the ladder polyominoes are Knutson. 
 We also show that when ladder polyominoes are ``based on" a parallelogram polyomino, the resulting polyomino is Knutson. 
 Additionally, we compute the reduced Gr\"{o}bner basis of these polyominoes in a special case.
 The techniques we develop to prove that ladder polyominoes are Knutson are used throughout this article. 
 Building upon this foundation, we prove that a special class of thin polyominoes is also Knutson. 
 Furthermore, we show that the set of all binomials associated with the inner intervals of these polyominoes forms the reduced Gr\"{o}bner basis for the polyomino ideals, and consequently, these polyominoes are in fact prime. Using results of \cite{JN24rookpolynomial} and \cite{MRR22:primality_of_polyomino}, we show that if every non-trivial intersection of two maximal inner intervals in a thin polyomino is a cell, then the polyomino is prime.
 Next, we show that under a certain condition, if a parallelogram polyomino is extracted from another parallelogram polyomino, the resulting collection of cells is Knutson, and we also compute the reduced Gr\"{o}bner basis. 

 The paper is organized as follows: in Section~\ref{sec:prelims}, we collect all necessary definitions and results on polyominoes and Knutson ideals. We begin Section~\ref{sec:thin1} by noting that, under certain conditions, ideals of K\"{o}nig type are indeed Knutson ideals. Using this observation, we show that polyomino ideals associated with simple thin, closed path, and weakly closed path polyominoes are Knutson ideals, by showing that these polyomino ideals are of K\"{o}nig type.
 In Section~\ref{sec:ladder}, we introduce the notion of ladder polyominoes, which are generalization of parallelogram polyominoes.
 Then, we show that ladder polyominoes are Knutson. 
 In Section~\ref{sec:thin2}, we prove that some other class of thin polyominoes are Knutson.
 Finally in Section~\ref{sec:collectioncells}, we inquire whether extracting a parallelogram polyomino from another parallelogram polyomino is Knutson or not. Under certain hypothesis, we prove that they are Knutson.
 
\section{Preliminaries}\label{sec:prelims}

\subsection{Polyominoes}
Let $(i,j),(k,l)\in \naturals^2$. We say that $(i,j)\leq(k,l)$ if $i\leq k$ and $j\leq l$. Consider $a=(i,j)$ and $b=(k,l)$ in $\naturals^2$ with $a\leq b$. The set $[a,b]=\{(m,n)\in \naturals^2: i\leq m\leq k,\ j\leq n\leq l \}$ is called an \textit{interval} of $\naturals^2$. 
In addition, if $i< k$ and $j<l$ then $[a,b]$ is a \textit{proper} interval. 
In such a case, we say $a, b$ the \textit{diagonal vertices} of $[a,b]$ and $c=(i,l)$, $d=(k,j)$ the \textit{anti-diagonal vertices} of $[a,b]$. 
Moreover, the elements $a$, $b$, $c$ and $d$ are called  the \textit{lower-left}, \textit{upper-right}, \textit{upper-left} and \textit{lower-right} \textit{vertices} of $[a,b]$ respectively.
 We also call $V([a,b])=\{a,b,c,d\}$ the set of \emph{vertices} of the interval, while the sets $\{a,c\},\{a,d\},\{b,d\},\{b,c\}$ are the \textit{edges} of the interval $[a,b]$.

An interval of the from $C = [a, a + (1, 1)]$ is called a \emph{cell} where $a\in \naturals$.
Let $\MP$ be a collection of cells. The \define{vertex set} of $\MP$ is $\bigcup_{C\in \MP} V(C)$.
Let $C$ and $D$ be two cells of $\MP$. Then, we say $C$ and $D$ are connected, if there is a sequence $\MC:C=C_1,\dots,C_m=D$ of cells (starting with $C$ and ending with $D$) of $\MP$ such that $C_i \cap C_{i+1}$ is an edge of both $C_i$ and $C_{i+1}$ for $i=1,\dots,m-1$. The sequence $\MC$ is called a \define{path} between $C$ and $D$ if $C_i \neq C_j$ for all $i\neq j$.
A \textit{polyomino} $\MP$ is a non-empty finite collection of cells in ${\naturals\setminus\{0\}}\times {\naturals\setminus\{0\}}$ where any two cells of $\MP$ are connected in $\MP$. 
A collection of cells (resp. polyomino) $\MP$ is called a \emph{sub-collection of cells} (resp. \emph{sub-polyomino}) of $\MP'$, if all the cells of $\MP$ are contained in $\MP'$.

 A collection of cells $\MP$ is called \textit{simple} if for any two cells $C$ and $D$ not in $\MP$ there exists a path of cells not in $\MP$ from $C$ to $D$. 
A polyomino $\MP$ is called \define{horizontally convex} (or row convex), if for any two cells $C, D$ of $\MP$ with lower-left
vertices $a = (i, j)$ and $b = (k, j)$ respectively, and such that $k > i$, it follows that
all cells with left lower corner $(l, j)$ with $i \leq l \leq k$ belong to $\MP$. Similarly, one
defines \define{vertically convex} (or column convex) polyominoes. The polyomino $\MP$ is called \define{convex} if it is horizontally
and vertically convex.
A convex polyomino $\MP$ is said to be a \emph{parallelogram polyomino} if $V(\MP)$ is a sub-lattice of $\naturals^2$.
A proper interval $[a, b]$, for $a,b\in V(\MP)$ is called an \emph{inner interval} of $\MP$ if all cells of $[a, b]$ belong to $\MP$.
An inner interval of $\MP$ is \emph{maximal} if it is maximal under inclusion.
 An interval $[a, b]$ with $a = (i, j)$ and $b = (k, l)$ with $k > i$ is called a \emph{horizontal edge interval} of $\MP$ if $j = l$ and the sets $\{(r, j),(r + 1, j)\}$ for $r = i, \ldots , k - 1$ are edges of cells of $\MP$. Similarly, one defines \emph{vertical edge interval} of $\MP$.

A collection of cells $\MP$ is \define{thin} if it does not have a $2 \times 2$ square such as the one shown in Figure~\ref{figure:twoRooks}. Let $C,D\in \MP$. We say that $C$ is a \define{neighbour} of $D$ if
$C\cap D$ is an edge of both $C$ and $D$. Let $\MP$ be a simple thin polyomino.
A cell $C$ is said to be an {\em end-cell} of a maximal inner interval $I$
if $C\in I$ and $C$ has exactly one neighbour cell in $I$.
A cell of $\MP$ is called {\em single} if it belongs to exactly one maximal
inner interval of $\MP$. 
Note that an inner interval of $\MP$ can take one of two forms: it either has vertices
$(i_1,j_1)$, $(i_1+1,j_1 )$, $(i_1,j_2 )$ and $(i_1+1, j_2 )$
or it has vertices
$(i_1,j_1)$, $(i_1,j_1+1)$, $(i_2,j_1 )$ and $(i_2, j_1+1)$
for some $i_1, i_2, j_1, j_2 \in \naturals$ with
$i_1 < i_2$ and $j_1 < j_2$.
By the \define{end-edges} of an inner interval, we refer to the edges $\{(i_1,j_1),(i_1+1,j_1)\}$ and $\{(i_1,j_2),(i_1+1, j_2)\}$ in the first case, and to the edges $\{(i_1,j_1),(i_1,j_1+1)\}$ and $\{(i_2,j_1),(i_2, j_1+1)\}$ in the second case.
\begin{figure}[h]
\begin{center}
\begin{tikzpicture}[scale=2]
\draw[] (0,0)--(0,1)--(1,1)--(1,0)--(0,0) (0,.5)--(1,.5) (.5,0)--(.5,1);
\end{tikzpicture}
\caption{}\label{figure:twoRooks}
\end{center}
\end{figure}

To a collection of cell $\MP$, Qureshi~\cite{QUR12} associated a finitely generated graded $K$-algebra $K[\MP]$ (over a field $K$) as follows:
Let $\MP$ be a collection of cells. Let $S = K[x_a \mid a\in V(\MP)]$ be the standard graded polynomial ring over a field $K$. 
For an inner interval $[a,b]$ of $\MP$, associate a binomial $x_{a} x_{b} - x_{c} x_{d},$ where $c$ and $d$ are the anti-diagonal vertices of $[a,b].$
Let $I_\MP$ be the binomial ideal of $S$ generated by the binomials $x_{a} x_{b} - x_{c} x_{d},$ where $[a,b]$ varies over all inner intervals of $\MP$.
We denote by $K[\MP]$ the quotient ring $S/I_\MP.$
When $\MP$ is a polyomino, the ideal $I_\MP$ is called the \emph{polyomino ideal} associated with $\MP$.
If $\MP$ is a simple polyomino, then it is known that $K[\MP ]$ is a Koszul Cohen-Macaulay normal domain~\cite[Corollary 2.2]{HM14},~\cite[Corollary~2.3]{QSS17}.

\subsection{Knutson ideals}

The notion of Knutson ideals was introduced by Conca and Varbaro after
Knutson’s work on compatibly split ideals~\cite{Knutson}. Our basic references on Knutson ideals are \cite{Knutson, Conca-Varbaro, Lisa1}.
\begin{definition}[\cite{Conca-Varbaro}]
	Let $K$ be a field and $S=K[x_1,\ldots,x_n]$ be a polynomial ring. Let $f\in S$ be a polynomial such that its leading term  $\ini_<f$ is a square-free monomial with respect to some monomial order $<$ on $S$. Define $\mathcal{C}_f$ to be the smallest set of ideals satisfying the following conditions:
	\begin{enumerate}
		\item $(f)\in \mathcal{C}_f$;
		\item If $I\in \mathcal{C}_f$ and $J$ is any ideal in $S$, then 
		$I : J \in \mathcal{C}_f$; hence  all the associated primes of
           $I$ are in $\mathcal{C}_f$.
		\item If $I$ and $J$ are ideals in $S$ such that  $I,J\in \mathcal{C}_f$ then  $I + J$ and $I\cap J$ must be in $\mathcal{C}_f$. 
	\end{enumerate}
	
	An ideal $I\subseteq S$ is called a \emph{Knutson ideal associated to $f$} if $I\in \mathcal{C}_f$. More generally,  an ideal $I$ is called a \emph{Knutson ideal} if $I\in \mathcal{C}_f$ for some $f$.
\end{definition}
 \begin{example}
 If $f=x_1\cdots x_n$ then the collection of Knutson ideals $\mathcal{C}_f$ is the set of all squarefree monomial ideals of $S$.
 \end{example}
  \begin{Remark}
      Since Knutson ideals are radical, one can replace the condition $(2)$ of the definition of Knutson ideals by the following: \\ 
      $(2)'$ If $I \in \mathcal{C}_f$ then 
      $\mathfrak{p}\in \mathcal{C}_f$ for every $\mathfrak{p}\in \Min(I)$, where $\Min(I)$ denotes the set of minimal primes of $I$.
      \end{Remark}
      \begin{lemma}
      \label{subset}
      	Let $K$ be a field and $S=K[x_1,\ldots,x_n]$ be a polynomial ring. Let $f,g\in S$ be  polynomials such that  $\ini_<(fg)$ is a square-free monomial with respect to some monomial order $<$ on $S$. Then $\mathcal{C}_f\subseteq \mathcal{C}_{fg}$.
      \end{lemma}

      \begin{proof}
      First note that since $\ini_<(fg)$ is a square-free monomial, 
      $\ini_<(f)$ is also a square-free monomial. We also have $(fg):(g)=(f)$, hence $(f)\in C_{fg}$. Hence from the construction of Knutson ideals, the proof follows.
      \end{proof}
When the characteristic of the field $K$ is prime and $\deg f=n$, then $f$ defines a splitting of the Frobenius homomorphism $F: S\to S$. Then one can see that all Knutson ideals associated to $f$ are compatibly split with respect to this map. Hence, Knutson ideals define $F$-split rings in prime characteristic. Moreover, by the assumption on initial term of $f$, Knutson ideals have squarefree initial ideals  in arbitrary characteristic, hence they are radical \cite{Knutson}.
As a consequence the extremal Betti numbers, Castelnuovo-Mumford regularity, and depth of Knutson ideals coincide with those of their
initial ideals. Furthermore, Gr\"obner bases of Knutson ideals are well-behaved with respect to sums as the union of Gr\"obner bases of Knutson ideals is a Gr\"obner basis of their sum. This fact makes computations of Gr\"obner bases of Knutson ideals easier in many cases. Thus study of Knutson ideals  become a useful tool to solve problems in applied algebra.

It has already been shown that some interesting classes of ideals are Knutson ideals for a suitable choice of $f$; for instance, ideals defining matrix Schubert varieties, ideals defining Kazhdan-Lusztig varieties \cite{Knutson}, determinantal ideals of Hankel matrices \cite{Lisa1} and generic matrices \cite{Lisa2}. Moreover building upon \cite{Lisa2} proved a graph is weakly closed if and only if its binomial edge ideal is a Knutson ideal \cite{Lisa-binomial}. Recently in \cite{DS-NB-S-V} it is shown that Ladder determinantal ideals are Knutson ideals.

\section{Knutson thin polyominoes I}\label{sec:thin1}

In this section, we start by proving that, under certain conditions, ideals of K\"{o}nig type are indeed Knutson ideals. We then apply this result to show that the polyomino ideals associated with simple thin, closed path, and weakly closed path polyominoes are Knutson ideals. We start by giving definition of ideals of K\"onig type.

\begin{definition}\cite[Definition\ 1.1]{HHM22konigtype}
Let $S=K[x_1,\ldots,x_n]$ be a polynomial ring  over a field $K$ and $I$ be a homogeneous ideal $S$ of height $h$. Then $I$ is of \emph{K\"{o}nig type} if it satisfies the following conditions:
\begin{enumerate}
 \item there exists a sequence of homogeneous polynomials $f_1,\ldots,f_h$ forming part of minimal system of homogeneous generators of $I$;
 \item there exists a monomial order $<$ on $S$ such that $\text{in}_<(f_1),\ldots,\text{in}_<(f_h)$ is a regular sequence.
\end{enumerate}
\end{definition}

We say that a polyomino $\MP$ is of  \emph{K\"{o}nig type} if the polyomino ideal $I_\MP$ is of K\"{o}nig type.

\begin{proposition}
\label{konig-implies-knutson}
Let $I$ be a homogeneous unmixed radical ideal of height $h$ in a polynomial ring $S$. If $I$ is of K\"{o}nig type with respect to $f_1,\ldots,f_h$, and $\ini_<(f_i)$ is square-free for all $i$, then $I$ is Knutson.
\end{proposition}
\begin{proof}
Let $f=f_1\cdots f_h$. Since $\ini_<(f)=\prod_{i=1}^{h}\ini_<(f_i)$, we get that $\ini_<(f)$ is square-free. 
Note that for each $1\leq i\leq h,$ $(f_i) = (f) : (f/f_i) \in C_f$.
Hence, $(f_1,\ldots,f_h)\in \mathcal{C}_f$. 
Since $I$ is unmixed, every minimal prime of $I$ is also a minimal prime of $(f_1,\ldots,f_h)$; hence, every minimal prime of $I$ is Knutson with respect to $f$. 
Since $I$ is radical, $I=\cap_{\mathfrak{p}\in \Min(I)} \mathfrak{p}$, where $\Min(I)$ denotes the set of minimal primes of $I$; thus $I$ is Knutson with respect to $f$.
\end{proof}

\subsection{Simple thin polyominoes}

Herzog and Hibi~\cite[Corollary\ 5.5]{HH23konig} showed that every simple thin polyomino is of K\"{o}nig type. 
Their argument relies on the repeated application of~\cite[Lemma\ 5.4]{HH23konig}. 
In the example provided below, we illustrate that \cite[Lemma\ 5.4]{HH23konig} in its original form is not sufficient to prove that every simple thin polyomino is of K\"{o}nig type. 
Instead, a more general version of~\cite[Lemma\ 5.4]{HH23konig} is required (see Lemma~\ref{lem:herzoghibi}), which must be used in a refined manner to complete the proof.

\begin{examplebox}\label{example:herzoghibi}

  \begin{figure}[h]%
		\centering
		\begin{tikzpicture}[scale=1]
		\draw[] 
        (1,5)--(1,3) (2,5)--(2,3)
        (3,2)--(3,4) (4,2)--(4,4)
        (5,3)--(5,4)

        (1,5)--(2,5)
        (1,4)--(5,4) (1,3)--(5,3)
        (3,2)--(4,2)
        ;
        \draw[red]
        (1,5)--(2,3)
        (1,3)--(3,4)
        (2,4)--(4,3)
        (3,3)--(4,4)

        ;

        \filldraw[] (1.5,4.5) circle (0pt) node[anchor=center]  {$C_1$};
        \filldraw[] (1.5,3.5) circle (0pt) node[anchor=center]  {$C_2$};
        \filldraw[] (2.5,3.5) circle (0pt) node[anchor=center]  {$C_3$};
        \filldraw[] (3.5,3.5) circle (0pt) node[anchor=center]  {$C_4$};
        \filldraw[] (4.5,3.5) circle (0pt) node[anchor=center]  {$C_5$};
        \filldraw[] (3.5,2.5) circle (0pt) node[anchor=center]  {$C_6$};

        \filldraw[] (1,5) circle (0.5pt) node[anchor=east]  {$a_1$};
        \filldraw[] (2,5) circle (0.5pt) node[anchor=west]  {$a_2$};
        \filldraw[] (1,4) circle (0.5pt) node[anchor=east]  {$a_3$};
        \filldraw[] (2,4) circle (0.5pt) node[anchor=south west]  {$a_4$};
        \filldraw[] (1,3) circle (0.5pt) node[anchor=east]  {$a_5$};
        \filldraw[] (2,3) circle (0.5pt) node[anchor=north]  {$a_6$};
        
        \filldraw[] (3,4) circle (0.5pt) node[anchor=south]  {$a_7$};
        \filldraw[] (3,3) circle (0.5pt) node[anchor=north east]  {$a_8$};
        \filldraw[] (4,4) circle (0.5pt) node[anchor=south]  {$a_9$};
        \filldraw[] (4,3) circle (0.5pt) node[anchor=north west]  {$a_{10}$};
        
        \filldraw[] (3,2) circle (0.5pt) node[anchor=north]  {$a_{11}$};
        \filldraw[] (4,2) circle (0.5pt) node[anchor=north]  {$a_{12}$};
        \filldraw[] (5,4) circle (0.5pt) node[anchor=west]  {$a_{14}$};
        \filldraw[] (5,3) circle (0.5pt) node[anchor=west]  {$a_{13}$};
		\end{tikzpicture}
		\caption{}%
	\label{fig:herzoghibi}
\end{figure}
Let $\MP$ be the polyomino as shown in Figure~\ref{fig:herzoghibi}. 
We aim to apply \cite[Lemma\ 5.4]{HH23konig} to examine whether 
$\MP$ is of K\"{o}nig type. 
The algorithm proposed in \cite{HH23konig} involves beginning at a cell of $\MP$ and 
then adding one cell at a time, with the goal of showing that the resulting polyomino is of K\"{o}nig type.
We start by considering the cell $C_1$. Let $\MP_1 = \{C_1\}$ and let $f = x_{a_3}x_{a_2} - x_{a_1}x_{a_4}$, with $<$ denoting the lexicographic order on $S_{\MP_1}$ as defined in \cite[Lemma\ 5.4]{HH23konig}, so that $\ini_<(f) = x_{a_1}x_{a_4}$. 
Thus, $\MP_1$ is of K\"{o}nig type with respect to $f$.

Define $\MP_i = \{C_1, \ldots, C_i\}$ for $2 \leq i \leq 4$. By \cite[Lemma\ 5.4]{HH23konig}, $\MP_i$ is of K\"{o}nig type for all $1 \leq i \leq 4$. Furthermore, $\MP_4$ is of K\"{o}nig type with respect to the generators 
$f_1 = x_{a_5}x_{a_2}-x_{a_1}x_{a_6},$
$f_2 = x_{a_5}x_{a_7}-x_{a_3}x_{a_8},$
$f_3 = x_{a_6}x_{a_9}-x_{a_4}x_{a_{10}},$
$f_4 = x_{a_8}x_{a_9}-x_{a_7}x_{a_{10}}.$ 
The initial terms of $f_i$'s are denoted by ``red" diagonal lines in Figure~\ref{fig:herzoghibi}.

We can extend $\MP_4$ to $\MP_5$ by adding either $C_5$ or $C_6$. In both cases, $\MP_5$ remains of K\"{o}nig type by~\cite[Lemma\ 5.4]{HH23konig}. 
Consider the case where $\MP_5 = \MP_4 \cup \{C_5\}$. In this case, $\MP_5$ is of K\"{o}nig type with respect to the generators $f_1,f_2, f_3, f' = x_{a_8}x_{a_{14}}-x_{a_7}x_{a_{13}},$ and $f'' = x_{a_{10}}x_{a_{14}}-x_{a_9}x_{a_{13}}.$
However, if we consider $\MP = \MP_5 \cup \{C_6\}$, we cannot use \cite[Lemma\ 5.4]{HH23konig} because $C_4$ does not satisfy the hypothesis of \cite[Lemma\ 5.4]{HH23konig}; specifically, $x_{a_8}x_{a_9} - x_{a_7}x_{a_{10}}$ is not part of the set of generators required to conclude that $\MP_5$ is of K\"{o}nig type. We encounter the same issue if we consider the case $\MP_5 = \MP_4 \cup \{C_6\}$. 
\end{examplebox}

\begin{Lemma}\label{lem:herzoghibi}
Let $\MP$ be a simple thin polyomino of K\"{o}nig type
with respect to $f_1, \ldots, f_r$ under a graded lexicographic order $<$ on $S$.
Let $I$ be a maximal inner interval of $\MP$ and $C$ be an end-cell of $I$, with $a, b\in V(C)$ such that $\{a,b\}$ is an edge of $C$ only in $\MP$ (see Figure~\ref{fig:herzoghibilemma}). 
Assume that $f_r=x_ax_d - x_bx_c$ and $\ini_< (f_r)=x_ax_d$ for some $c,d\in V(I)$. 
Let $\MP'$ be a simple thin polyomino obtained by adding a cell $D$ to $\MP$  with $V(C)\cap V(D) = \{a,b\}$. 
Write $V(D) = \{a,b, u,v\}$ such that $\{a,u\}$ forms an edge of $D$.
Then $\MP'$ is of K\"{o}nig type with respect to $f_1, \ldots, f_{r-1}, f', f'',$ where $f' = x_ux_d - x_vx_c , f'' = x_ux_b - x_vx_a$, under a graded lexicographic order $<'$ on $S'=S[x_u,x_v]$ with $\ini_{<'} (f')=x_ux_d$, $\ini_{<'} (f'')=x_ax_v$.
\end{Lemma}

 \begin{figure}[h]%
		\centering
		\begin{tikzpicture}[scale=1]
		\draw[] 
        (1,5)--(1,2) (2,5)--(2,2)  
        (-2,4)--(-2,3) (-1,4)--(-1,3)
        (-0.25,4)--(-0.25,3)

        (1,5)--(2,5)
        (1,2)--(2,2)
        (0.5,4)--(2,4) (0.5,3)--(2,3)
        (-2,4)--(-1,4) (-2,3)--(-1,3)
        
        ;
        \draw[red]
        (-0.25,3)--(2,4)
        ;
        \draw[dotted, thick]
        (2,3)--(3,3)--(3,4)--(2,4)
        (0.5,4)--(-1,4) (0.5,3)--(-1,3)

        ;

        \filldraw[] (1.5,3.5) circle (0pt) node[anchor=center]  {$C$};
        \filldraw[] (2.5,3.5) circle (0pt) node[anchor=center]  {$D$};
        \filldraw[] (-2.5,3.5) circle (0pt) node[anchor=center]  {$I$};

        \filldraw[] (2,4) circle (0.5pt) node[anchor=south west]  {$a$};
        \filldraw[] (2,3) circle (0.5pt) node[anchor=north west]  {$b$};
        
        \filldraw[] (-0.25,4) circle (0.5pt) node[anchor=south]  {$c$};
        \filldraw[] (-0.25,3) circle (0.5pt) node[anchor=north]  {$d$};

        \filldraw[] (3,4) circle (0.5pt) node[anchor=west]  {$u$};
        \filldraw[] (3,3) circle (0.5pt) node[anchor=west]  {$v$};
		\end{tikzpicture}
		\caption{}%
	\label{fig:herzoghibilemma}
\end{figure}

\begin{proof}
    The proof follows the argument of~\cite[Lemma\ 3.2, 3.3]{HH23konig}. 
\end{proof}

Before stating our main theorem of this section, we recall the notion of collapsing a simple thin polyomino in a maximal inner interval, as introduced by Rinaldo and Romeo~\cite[Proposition 3.7]{RR21}. We define this notion as it appears in \cite[Definition\ 3.2]{KV22}.

\begin{definition}
\label{def:coll}
Let $\MP$ be a simple thin polyomino. A \define{collapse datum} on $\MP$
is a tuple $(I, J,\MP^I)$, where $I$ and $J$ are maximal inner
intervals and $\MP^I$ is a sub-polyomino of $\MP$
satisfying the following conditions:
\begin{enumerate}

  \item
$J$ is the only maximal inner interval of $\MP$ such that $I \cap J$ is a cell;

\item
$\MP^I \subseteq J$ and
$I \cap J \not \subset \MP^I$.

\item
$\MP \setminus ( I \cup \MP^I)$ is a non-empty sub-polyomino of $\MP$.
\end{enumerate}
\end{definition}

When $\MP$ has at least two maximal inner intervals, the maximal inner intervals $I$ and $J$ defined in the Definition~\ref{def:coll} exist by \cite[Lemma\ 3.6]{RR21}.

\begin{theorem}\label{thm:simplethin}
    Every simple thin polyomino is of K\"{o}nig type. 
\end{theorem}

\begin{proof}
    Let $\MP$ be a simple thin polyomino with $m$ number of cells.
    We proceed by induction on the number of maximal inner intervals of $\MP$ and show the following: 
     $(\#)$ $\MP$ is of K\"{o}nig type with respect to $f_1, \ldots, f_m$ under a graded lexicographic order $<$ on $S$ such that
       \begin{asparaenum}
       \item
           if a cell $A$ is both a single cell and an end-cell of a maximal inner interval $I$ of $\MP$, then there exists an $i\in [m]$ such that $f_i=I_A$, where $I_A$ is the polyomino ideal of the cell $A$;
       
       \item
            for every maximal inner interval $I$ of $\MP$ and every end-edge $e$ of $I$, either there exists a vertex $v\in e$ such that $\gcd(x_v,\ini_<(f_i)) =1$ for all $1\leq i\leq m$ or there exists an $i\in [m]$ such that $f_i$ correspond to an inner interval in $I$ containing edge $e$.
       \end{asparaenum}

   \begin{figure}[h]%
		\centering
		\begin{tikzpicture}[scale=1]
		\draw[] 
        (1,0)--(2,0)--(2,1)--(1,1)--(1,0)
        (4,0)--(5,0)--(5,1)--(4,1)--(4,0)
        ;
        \draw[dotted, thick] (2,0)--(4,0) (2,1)--(4,1)
        (3,1)--(3,0); 
        
        \filldraw[] (4.5,.5) circle (0pt) node[anchor=center]  {$A_m$};
        \filldraw[] (1.5,0.5) circle (0pt) node[anchor=center]  {$A_1$};
        
        \filldraw[] (5,1) circle (0.5pt) node[anchor=west]  {$u_{m+1}$};
        \filldraw[] (5,0) circle (0.5pt) node[anchor=west]  {$v_{m+1}$};
        \filldraw[] (3,0) circle (0.5pt) node[anchor=north]  {$v_i$};
        \filldraw[] (3,1) circle (0.5pt) node[anchor=south]  {$u_i$};
        
        \filldraw[] (1,0) circle (0.5pt) node[anchor=east]  {$v_1$};
        \filldraw[] (1,1) circle (0.5pt) node[anchor=east]  {$u_1$};

		\end{tikzpicture}
		\caption{}%
	\label{fig:interval}
\end{figure}

   First, assume that $\MP$ has exactly one maximal inner interval. Then, $\MP$ is as shown in Figure~\ref{fig:interval}. 
   Under the notations of Figure~\ref{fig:interval}, $S = K[x_{v_1},\ldots,x_{v_{m+1}},x_{u_1},\ldots,x_{u_{m+1}}]$.
   Let $<$ be the lexicographic order on $S$ induced by the ordering $x_{v_1}< \cdots < x_{v_{m+1}}< x_{u_1}< \cdots < x_{u_{m+1}}$.
   For $i\in [m]$, 
   define $f_i =  x_{v_{i}}x_{u_{i+1}} - x_{u_{i}}x_{v_{i+1}}$; then, $\ini_< (f_i) = x_{u_{i+1}}x_{v_{i}}$.
   Clearly, $\MP$ is of K\"{o}nig type with respect to $f_1, \ldots, f_m$ and satisfies $\#$.
   
    So we may assume that $\MP$ has at least two maximal inner interval. 
    Let $(I,J,\MP^I)$ be a collapse datum of $\MP$ and let $I\cap J = \{C\}$.
    Let $C_1,\ldots, C_r, D_1,\ldots,D_s$ be the cells of $I$ different from $C$. 
    Suppose that $C_1$ and $D_1$ are neighbors of $C$. 
    Furthermore, for all $i \in {1,\ldots, r-1}$, $C_i$ and $C_{i+1}$ are neighbours, and for all $j \in {1,\ldots, s-1}$, $D_j$ and $D_{j+1}$ are neighbours.
    Note that either $\{C_1,\ldots, C_r\}$ or $\{D_1,\ldots,D_s\}$ may be empty;  
    without loss of generality, we may assume that $\{C_1,\ldots, C_r\}$ is non-empty.
    In the case where $\MP^I\neq \varnothing$, let $F_1,\ldots, F_t$ be the cells of $\MP^I$. 
    Assume that $C$ and $F_1$ are neighbours and for all $i\in \{1,\ldots, t-1\}$, $F_i$ and $F_{i+1}$ are neighbours.
    Moreover, for all $1\leq i\leq t$, assume that $V(F_i) = \{u_i,v_i, u_{i+1}, v_{i+1}\}$ such that $V(C)\cap V(F_1) = \{v_1,u_1\},$ $V(F_i)\cap V(F_{i+1}) = \{v_{i+1},u_{i+1}\}$ and $\{u_{i},u_{i+1}\}$ is an edge of $F_i$ for all $1\leq i\leq t$.
    See Figure~\ref{fig:simplethin} for a schematic showing of the cells $C, {C_i}'s, {D_j}'s$ and ${F_k}'s$. 

    \begin{figure}[h]
		\centering
		\begin{tikzpicture}[scale=.9]
		\draw[] 
        (1,0)--(4,0)--(4,1)--(1,1)--(1,0)
        (2,0)--(2,2)--(3,2)--(3,0)
        (2,2)--(3,2)
        (5,0)--(6,0)--(6,1)--(5,1)--(5,0)
        (3,3)--(3,4)--(2,4)--(2,3)--(3,3)
        (3,0)--(3,-1)--(2,-1)--(2,0)
        (3,-2)--(3,-3)--(2,-3)--(2,-2)--(3,-2)
        ;
        \draw[dotted, thick] (4,0)--(5,0) (4,1)--(5,1)
        (3,2)--(3,3) (2,2)--(2,3)
        (3,-1)--(3,-2) (2,-1)--(2,-2); 
        
        \draw[dashed] (2,1)--(1,1.5);
        \draw[dashed] (2,0)--(1,-.5);
        \draw[dashed] (1,1.5) arc(90:270:1);
        \draw[mybrace=0.5, thick] (3,4.1) -- (6,4.1);
		
        \filldraw[] (2.5,4.5) circle (0pt) node[anchor=center]  {$I$};
        \filldraw[] (4.5,4.5) circle (0pt) node[anchor=center]  {$\MP^I$};
        \filldraw[] (2.5,.5) circle (0pt) node[anchor=center]  {$C$};
        
        \filldraw[] (2.5,3.5) circle (0pt) node[anchor=center]  {$C_r$};
        \filldraw[] (2.5,1.5) circle (0pt) node[anchor=center]  {$C_1$};
        \filldraw[] (2.5,-.5) circle (0pt) node[anchor=center]  {$D_1$};
        \filldraw[] (2.5,-2.5) circle (0pt) node[anchor=center]  {$D_s$};
        
        \filldraw[] (3.5,.5) circle (0pt) node[anchor=center]  {$F_1$};
        \filldraw[] (5.5,.5) circle (0pt) node[anchor=center]  {$F_t$};

        \filldraw[] (6,1) circle (0.5pt) node[anchor=west]  {$u_{t+1}$};
        \filldraw[] (6,0) circle (0.5pt) node[anchor=west]  {$v_{t+1}$};
        \filldraw[] (4.5,0) circle (0.5pt) node[anchor=north]  {$v_i$};
        \filldraw[] (4.5,1) circle (0.5pt) node[anchor=south]  {$u_i$};
        
        \filldraw[] (3,0) circle (0.5pt) node[anchor=north west]  {$v_1$};
        \filldraw[] (3,1) circle (0.5pt) node[anchor=south west]  {$u_1$};
        
		\end{tikzpicture}
		\caption{}\label{fig:simplethin}
\end{figure}

    Consider the sub-polyomino $\MP' = (\MP \setminus (I \cup \MP^I))\cup \{C\}$ of $\MP$. 
    Note that the number of maximal inner intervals of $\MP'$ are exactly one less than the number of maximal inner intervals of $\MP$.
    So, by induction hypothesis, $\MP'$ is of K\"{o}nig type and satisfies $\#$.
   
    Note that $C$ is both a single cell and an end-cell of the maximal interval containing it in $\MP'$, which is $J\setminus \MP^I$.
    Since $\MP'$ satisfies $(1)$ of $\#$, the polyomino $\MP'\cup \{C_1\}$ is of K\"{o}nig type and satisfies $\#$ by Lemma~\ref{lem:herzoghibi}.
    Repeatedly applying Lemma~\ref{lem:herzoghibi}, we deduce that $\MP'\cup \{C_1,\ldots,C_r\}$ is of K\"{o}nig type and satisfies $\#$.

    Consider the case when $\MP^I\neq \varnothing$. 
    Since $C$ is non-single end-cell of the maximal inner interval $J\setminus \MP^I$ of $\MP'\cup \{C_1,\ldots,C_r\}$, $(2)$ of $\#$ holds for $\{v_1,u_1\}$. 
    First suppose that, the end-edge $\{v_1,u_1\}$ of $J\setminus \MP^I$ contains a vertex, say $v_1$, such that $\gcd(x_{v_1},\ini_<(f_i)) =1$ for all $1\leq i\leq m$.
    Extend the graded lexicographic order $<$ on $S$ to $S_1 = S[x_{v_2},\ldots,x_{v_{m+1}},x_{u_2},\ldots,x_{u_{m+1}}]$ by the total order $x_{v}< x_{v_2}< \cdots < x_{v_{t+1}}< x_{u_2}< \cdots < x_{u_{t+1}}$, for all $v\in \MP'\cup \{C_1,\ldots,C_r\}$. 
    For all $1\leq i \leq t$, define $f_{m+i} =x_{u_{i+1}}x_{v_i}-x_{u_{i}}x_{v_{i+1}}$;
    then, $\ini_<(f_{m+i})=x_{u_{i+1}}x_{v_i}$ for all $1\leq i \leq t$. Thus, $\MP'\cup \{C_1,\ldots,C_r,F_1,\ldots,F_t\}$ is of K\"{o}nig type with respect to $f_1, \ldots, f_{m+t}$.
   Next, suppose that there exists an $i\in [m]$ such that $f_i$ correspond to an inner interval in $J\setminus \MP^I$ containing edge $\{v_1,u_1\}$. Then, the polyomino $\MP'\cup \{C_1,\ldots,C_r,F_1\}$ is of K\"{o}nig type and satisfies $\#$ by Lemma~\ref{lem:herzoghibi}. Again, by repeatedly using Lemma~\ref{lem:herzoghibi}, we get that $\MP'\cup \{C_1,\ldots,C_r,F_1, \ldots, F_s\}$ is of K\"{o}nig type and satisfies $\#$.

    Consider the case when $\{D_1,\ldots,D_s\}\neq \varnothing$.
    Since $\MP'\cup \{C,C_1,\ldots,C_r,F_1, \ldots, F_s\}$ satisfies $(2)$ of $\#$, the polyomino $\MP'\cup \{C,C_1,\ldots,C_r,F_1, \ldots, F_s,D_1\}$ is of K\"{o}nig type and satisfies $\#$ by Lemma~\ref{lem:herzoghibi}. Proceeding similarly, we get that $\MP'\cup \{C,C_1,\ldots,C_r,D_1, \ldots, D_s,F_1,\ldots,F_t\}$ is of K\"{o}nig type and satisfies $\#$.
    Hence the proof.
\end{proof}

We now illustrate above theorem for the polyomino as shown in Figure~\ref{fig:herzoghibi}.
Using proof of Theorem~\ref{thm:simplethin}, we get $\MP$ is of K\"{o}nig type with respect to the generators 
$f_1 = x_{a_3}x_{a_2}-x_{a_1}x_{a_4},$
$f_2 = x_{a_5}x_{a_7}-x_{a_3}x_{a_8},$
$f_3 = x_{a_{6}}x_{a_{14}}-x_{a_4}x_{a_{13}},$ 
$f_4 = x_{a_{11}}x_{a_9}-x_{a_7}x_{a_{12}},$ 
$f_5 = x_{a_{11}}x_{a_{10}}-x_{a_8}x_{a_{12}}$ and 
$f_6 = x_{a_{10}}x_{a_{14}}-x_{a_9}x_{a_{13}}$. 
Moreover, $\ini_<(f_1) = x_{a_{2}}x_{a_{3}}, \ini_<(f_2) = x_{a_{5}}x_{a_{7}}, \ini_<(f_3) = x_{a_{4}}x_{a_{13}}, 
\ini_<(f_4) = x_{a_{9}}x_{a_{11}}, \ini_<(f_5) = x_{a_{8}}x_{a_{12}}$ and $\ini_<(f_6) = x_{a_{10}}x_{a_{14}}$.

\begin{corollary}
    Let $\MP$ be a simple thin polyomino. Then, $I_\MP$ is Knutson.
\end{corollary}
\begin{proof}
    By Theorem~\ref{thm:simplethin}, there exist generators $f_1,\ldots,f_m$ of $I_\MP$ corresponding to distinct inner intervals of $\MP$ such that $I_\MP$ is of K\"{o}nig type with respect to $f_1,\ldots,f_m$. Take $f = f_1\cdots f_m.$ Then, $I_\MP\in C_f$ by Proposition~\ref{konig-implies-knutson}.
\end{proof}

\subsection{Closed path and weakly closed path polyominoes}

We begin by recalling the definitions of closed path polyominoes and weakly closed path polyominoes as defined in \cite{CN23closedpathprimality} and \cite{CNU22weaklyclosedprimality} respectively. We then use~\cite{DN22konigtype} to show that the polyomino ideal of a closed path polyomino is Knutson. Subsequently, we prove that the polyomino ideal of a weakly closed path polyomino is of K\"{o}nig type. Recently, this result was also proved by Navarra~\cite{navarra2024weaklyclosedkonig}; our argument is also on the same line of~\cite{navarra2024weaklyclosedkonig}. Finally, we prove that polyomino ideal of a weakly closed path polyomino is Knutson.

\begin{definition}
Let $\mathcal{P}$ be a finite non-empty collection of cells. Then $\mathcal{P}$ is called a \emph{closed path polyomino} if it is a sequence of cells $A_0,\ldots,A_{n}$ with $n\geq 6$ satisfying the following properties:
\begin{enumerate}
\item $A_0 = A_{n}$;
\item $A_i\neq A_j$ , for all $i \neq j$ and  $i, j \in  \{1, \cdots, n\}$;
\item $A_{i}\cap A_{i+1}$ is a common edge for all $i\in \{0, \ldots,n-1\}$;
\item  $V(A_i)\cap V(A_j)=\emptyset$ for all $i\in \{0,\ldots,n-1\}$ and for all $j\notin \{i-2,i-1,i,i+1,i+2\}$, where the indices are reduced modulo $n$.
\end{enumerate}
\end{definition}

\begin{definition}
 Let $\mathcal{P}$ be a finite non-empty collection of cells. Then $\mathcal{P}$ is called \emph{weakly closed path} if it is a sequence of cells $A_0, A_1,\ldots,A_{n-1}, A_n$ with $n\geq 6$ satisfying the following properties:
 \begin{enumerate}
 \item $A_0 = A_n$;
  \item $A_i\neq A_j$ , for all $i \neq j$ and  $i, j \in  \{1, \cdots, n\}$;
  \item $A_0\cap A_{n-1}$ is a vertex;
  \item $A_{i}\cap A_{i+1}$ is a common edge for all $i\in \{0, \ldots,n-2\}$;
  \item $V(A_i)\cap V(A_j)=\emptyset$ for all $i\in \{1,\ldots,n\}$ and for all $j\notin \{i-2,i-1,i,i+1,i+2\}$, where the indices are reduced modulo $n$.
 \end{enumerate}
\end{definition}
\begin{Remark}
    In a closed path, one can order the cells such that every cell has an edge in common with its consecutive cell. In a weakly closed path the same holds, with the exception of exactly two consecutive cells that have just a vertex in common. Note that each closed path polyomino and weakly closed path polyominoe has a unique hole.
\end{Remark}

\begin{proposition}
    Let $\mathcal{P}: A_0,\ldots,A_{n-1}$ be a closed path polyomino.
 Then $I_{\mathcal{P}}$ is Knutson. 
\end{proposition}

\begin{proof}
    By~\cite[Theorem\ 4.9]{DN22konigtype}, $I_{\mathcal{P}}$ is of K\"{o}nig type.
   To prove this, the authors introduced an total order $<$ on the vertices of $\mathcal{P}$.
   For each cell $A_i$ of $\mathcal{P}$, they assigned a generator $f_i$ of $I_\MP$ corresponding to an inner interval of $\MP$ such that $\gcd(\ini_<(f_i), \ini_<(f_j))=1$ for all $i\neq j$.
   Consider $f=\prod f_i$, then $\ini_<(f)$ is square-free. 
   Note that $I_{\mathcal{P}}$ is unmixed by \cite[Theorem\ 4.19]{CNV2024polyocollection}, 
   and by~\cite{CNU22grobnerclosedpath}, there exists a Gr\"{o}bner basis of $I_{\mathcal{P}}$ with square-free initial terms, so $I_{\mathcal{P}}$ is radical.  
   Hence the proof follows from Proposition~\ref{konig-implies-knutson}.
\end{proof}

\begin{theorem}\label{thm:weaklyclosedpathkonig}
 Let $\mathcal{P}: A_0,\ldots,A_{n-1}$ be a weakly closed path polyomino.
 Then $I_{\mathcal{P}}$ is of K\"{o}nig type.
 \end{theorem}

\textbf{Idea of the proof:}
Our proof follows similar strategies of \cite{DN22konigtype}.
Let $\mathcal{P}$ be a weakly closed path polyomino with $n$ distinct cells $A_0, A_1,\ldots, A_{n-1}$. Since $\mathcal{P}\setminus\{A_0\}$ is simple, by \cite[Proposition\ 4.10]{CNV2024polyocollection}, one can see that $\height(I_{\mathcal{P}})=\rank{(\mathcal{P})}= n$. In order to show that $I_{\mathcal{P}}$ is of K\"{o}nig type, first, we give a suitable ordering on the vertices of $\mathcal{P}$, then we will consider the lexicographic ordering $<_{lex}$ induced by this total order on $\{x_v : v \in V (\mathcal{P})\}$. Next, we choose $n$ generators $f_1, \ldots, f_n$ of $I_{\mathcal{P}}$ such that their initial terms do not share any common variables. Hence  $\ini_{<_{lex}}(f_1), \ldots , \ini_{<_{lex}} (f_n)$ forms a regular sequence. Therefore, $I_{\mathcal{P}}$ will be of K\"onig type.

\begin{table}       
\tiny
         \begin{center}
\begin{tabular}{ |c|c|c|} 

&
IF it occurs ... 
&
THEN we refer to ...\\
  \hline
  I
  &
  \begin{tikzpicture}[scale=.7]
   \draw[](0,0)--(3,0)--(3,1)--(0,1)--(0,0);
   \draw[](1,0)--(1,1);
    \draw[](2,0)--(2,1);
    \draw[blue](0,1)--(1,0);
    \filldraw[color=blue] (0,1) circle (2pt);
     \filldraw[color=blue] (1,0) circle (2pt);
     \node at (0.5,0.5) {$A_{i}$};
      \node at (1.5,0.5) {$A_{i+1}$};
      \node at (2.5,0.5) {$A_{i+2}$};
      \filldraw[] (1,0) circle (0pt) node[anchor=north]{$t'$};
            
    \filldraw[] (0,1) circle (0pt) node[anchor=south]{$t$};
     \end{tikzpicture} 
     
     &
     
     \begin{tikzpicture}[scale=.7]
   \draw[](0,0)--(3,0)--(3,1)--(0,1)--(0,0);
   \draw[](1,0)--(1,1);
    \draw[](2,0)--(2,1);
    \draw[blue](0,1)--(1,0);
    \draw[orange](1,1)--(2,0);
    \filldraw[color=blue] (0,1) circle (2pt);
     \filldraw[color=blue] (1,0) circle (2pt);
     \filldraw[color=orange] (1,1) circle (2pt);
     \filldraw[color=orange] (2,0) circle (2pt);
      \node at (0.5,0.5) {$A_{i}$};
      \node at (1.5,0.5) {$A_{i+1}$};
      \node at (2.5,0.5) {$A_{i+2}$};

      \filldraw[] (1,0) circle (0pt) node[anchor=north]{$t'$};
      \filldraw[] (2,0) circle (0pt) node[anchor=north]{$(t+1)'$};
            
    \filldraw[] (0,1) circle (0pt) node[anchor=south]{$t$};
    \filldraw[] (1,1) circle (0pt) node[anchor=south]{$t+1$};
         \end{tikzpicture}\\
\hline  
\end{tabular}
\caption{}
\label{table}
\end{center}
\end{table}

In order to define a suitable total ordering
on $\{x_v : v \in V (\mathcal{P})\}$, we order the vertices inductively using Table~\ref{table} (Table of \cite{DN22konigtype}) below. We start with the cell $A_0$.  Let $Y^{(1)}=Y^{(1)}_1\sqcup Y^{(1)}_2$, where $Y^{(1)}_1=\{1\}, Y^{(1)}_2=\{1'\}$.  Let $j\geq 2$ and assume that $Y^{(j-1)}$ is known. We want to define $Y^{(j)}$.
We refer to Table~\ref{table} up to suitable rotations and reflections of 
$\mathcal{P}$. If one of the configurations in the left column of Table~\ref{table} occurs, where the blue vertices are in $Y^{(j-1)}$,
then we denote by $k$ the maximum integer such that $t + k$ 
is a red vertex in the picture displayed
in the corresponding right column. Then, we set $Z^{(j)}_1 = 
\{x_{t+1}, \cdots, x_{t+k}\}$ and $Z^{(j)}_2 = \{x_{(t+1)'}, \cdots, x_{(t+k)'}\}$.
For all $a \in Y^{(j-1)}_1$, we put $x_{1'}<_{lex} x_{h_2}<_{lex} x_{h_1}<_{lex} x_a$
with $t < h_1 <h_2 \leq t + k$, and for all $b\in Y_2^{(j-1)},$
we put $x_{m_2}<_{lex}x_{m_1}<_{lex}x_b$
with $t'<m_1<m_2\leq (t+k)'$.
Therefore,  we define $Y^{(j)} = Y^{(j)}_1 \sqcup Y^{(j)}_2,$
where
$Y^{(j)}_1 = Y^{(j-1)}_1 \sqcup Z^{(j)}_1$
and $Y^{(j)}_2 = Y^{(j-1)}_2 \sqcup Z^{(j)}_2$. 
Finally, we argue that this process allows us to complete the labeling of all the vertices of $\mathcal{P}$.

\begin{proof}[Proof of Theorem~\ref{thm:weaklyclosedpathkonig}]
Note that in the neighbourhood of $A_0$, $\mathcal{P}$ is one of the following types in Figure~\ref{fig:weakly-closed-1}.

\begin{figure}[h]
\small
    \begin{subfigure}[t]{4.5cm}
    \centering
		\begin{tikzpicture}[scale=.8]
  \draw[](0,0)--(2,0)--(2,1)--(0,1)--(0,0);
   \draw[](1,0)--(1,1);
    \node at (0.5,0.5) {$A_0$};
    \node at (1.5,0.5) {$A_1$};
    \draw[](0,0)--(-1,0)--(-1,-1)--(0,-1)--(0,0);
    \node at (-0.5,-0.5) {$A_{n-1}$};
    
    \node at (-1.5,-0.5) {$\cdots$};
    \node at (-0.5,-1.5) {$\vdots$};
    
     \node at (2.5,0.5) {$\cdots$};
    \node at (1.5,-0.5) {$\vdots$};
    \node at (1.5, 1.5) {$\vdots$};
	\end{tikzpicture}
    \caption{}%
    \end{subfigure}
    \quad
     \begin{subfigure}[t]{4.5cm}
		\centering
       \begin{tikzpicture}[scale=.8]
           \draw[](0,0)--(0,2)--(1,2)--(1,0)--(0,0);
            \draw[](0,1)--(1,1);
             \draw[](0,0)--(-1,0)--(-1,-1)--(0,-1)--(0,0);
             \node at (0.5,0.5) {$A_0$};
              \node at (0.5,1.5) {$A_1$};
               \node at (-0.5,-0.5) {$A_{n-1}$};
               
                \node at (-1.5,-0.5) {$\cdots$};
                 \node at (-0.5,-1.5) {$\vdots$};
                 
                 \node at (1.5,1.5) {$\cdots$};
               \node at (0.5,2.5) {$\vdots$};
               \node at (-0.5, 1.5) {$\cdots$};
               
         \end{tikzpicture}
        \caption{}
	\end{subfigure}
    \caption{}%
    \label{fig:weakly-closed-1}
\end{figure}

  Depending on the position of $A_1$ relative to $A_0$, we have two different ways to label the vertices of $A_0$. When $A_1$ is positioned to the right of $A_0$, we label vertices of $A_0$ as in Figure~\ref{fig:weakly-closed-2(a)}. When $A_1$ is positioned above $A_0$, we label vertices of $A_0$ as in Figure~\ref{fig:weakly-closed-2(b)}. We set  $Y^{(1)}=Y_1^{(1)}\sqcup Y_2^{(1)}$, where $Y_1^{(1)}=\{x_1\}$ and $Y_2^{(1)} = \{x_{1'}\}$.
 
    \begin{figure}[h]
    \small
    \begin{subfigure}[t]{4.5cm}
    \centering
		\begin{tikzpicture}[scale=.8]
  \draw[](0,0)--(2,0)--(2,1)--(0,1)--(0,0);
   \draw[](1,0)--(1,1);
   
    \node at (0.5,0.5) {$A_0$};
    \node at (1.5,0.5) {$A_1$};
    \draw[](0,0)--(-1,0)--(-1,-1)--(0,-1)--(0,0);
    \node at (-0.5,-0.5) {$A_{n-1}$};
    
    \node at (-1.5,-0.5) {$\cdots$};
    \node at (-0.5,-1.5) {$\vdots$};
    
     \node at (2.5,0.5) {$\cdots$};
    \node at (1.5,-0.5) {$\vdots$};
    \node at (1.5, 1.5) {$\vdots$};

    \draw[](1,0)--(0,1);
    \node at (-0.3,1.2) {$1$};
    \node at (0.9,-0.3) {$1'$};

     \filldraw[color=blue] (0,1) circle (2pt);
     \filldraw[color=blue] (1,0) circle (2pt);
	\end{tikzpicture}
    \caption{}
     \label{fig:weakly-closed-2(a)}%
    \end{subfigure}
    \quad
     \begin{subfigure}[t]{4.5cm}
		\centering
       \begin{tikzpicture}[scale=.8]
           \draw[](0,0)--(0,2)--(1,2)--(1,0)--(0,0);
            \draw[](0,1)--(1,1);
             \draw[](0,0)--(-1,0)--(-1,-1)--(0,-1)--(0,0);
             \node at (0.5,0.5) {$A_0$};
              \node at (0.5,1.5) {$A_1$};
               \node at (-0.5,-0.5) {$A_{n-1}$};
               
                \node at (-1.5,-0.5) {$\cdots$};
                 \node at (-0.5,-1.5) {$\vdots$};
                 
                 \node at (1.5,1.5) {$\cdots$};
               \node at (0.5,2.5) {$\vdots$};
               \node at (-0.5, 1.5) {$\cdots$};

            \filldraw[color=blue] (0,1) circle (2pt);
             \filldraw[color=blue] (1,0) circle (2pt);

              \draw[](1,0)--(0,1);
             
             \node at (-0.3,1.2) {$1'$};
             \node at (0.9,-0.3) {$1$};
         \end{tikzpicture}
        \caption{}\label{fig:weakly-closed-2(b)}
	\end{subfigure}
    \caption{}%
\end{figure}

 We then proceed with labeling the vertices according to the Table~\ref{table}, up to rotations and reflections, until $A_{n-3}$ or $A_{n-2}$, depending on the configuration of the cells $A_{n-4}, A_{n-3},A_{n-2}, A_{n-1}$, and the labelling of $A_{n-3}$. Note that we have two possibilities of $A_{n-1}$ and $A_{n-2}$, i.e., $A_{n-2}$ could be either to the left of $A_{n-1}$ or below $A_{n-1}$.

The blue vertices denote the vertices which have been labelled in the penultimate step and red vertices denote the vertices that are labelled in this step to complete the labeling.

First, we analyze the case where $A_{n-2}$  is positioned to the left of $A_{n-1}$. If $A_{n-3}$ is also on the left to $A_{n-2}$, we proceed to complete the labeling of the vertices of $\mathcal{P}$ as shown in Figure~\ref{fig:weakly-closed-5}. In this figure, blue vertices represent those that were labeled in the previous step, while red vertices denote those labeled in the current step to complete the labeling.

\begin{figure}[h]
\scriptsize
    \begin{subfigure}[t]{8cm}
    \centering
 
      \caption{}
	\end{subfigure}
    \quad
     \begin{subfigure}[t]{7cm}
		\centering
%
 
 \caption{}
	\end{subfigure}
    \quad
     \begin{subfigure}[t]{7cm}
		\centering
     %
      \caption{}
	\end{subfigure}
    \quad
     \begin{subfigure}[t]{7cm}
		\centering
%
 
      \caption{}
	\end{subfigure}
       \caption{}
     \label{fig:weakly-closed-5}
 \end{figure} 
 
     Next, we consider the subcase when $A_{n-3}$ is located below of $A_{n-2}$. In this scenario, 
     Figure~\ref{fig:weakly-closed-6} represents the vertex labeling for all possible configurations that may arise.
     
     \begin{figure}[h]
     \scriptsize
    \begin{subfigure}[t]{8cm}
    \centering
%
 
\caption{}
	\end{subfigure}
       \caption{}
     \label{fig:weakly-closed-6}
 \end{figure} 

For the subcase where $A_{n-3}$ is positioned above of $A_{n-2}$. Figure~\ref{fig:weakly-closed-7} presents all possible configurations.

\begin{figure}[h]
\scriptsize
    \begin{subfigure}[t]{8cm}
    \centering
%
 
\caption{}
	\end{subfigure}
       \caption{}
     \label{fig:weakly-closed-7}
 \end{figure}

Similarly one can conclude the second case, i.e., $A_{n-2}$ is positioned below $A_{n-1}$. All possible configurations for this case are shown in Figure~\ref{fig:weakly-closed-8}.

\begin{figure}[h]
\scriptsize
      \begin{subfigure}[t]{5cm}
    \centering
%
   \caption{}
	\end{subfigure}
     \quad
     \begin{subfigure}[t]{5cm}
		\centering
     %
             \caption{}
	\end{subfigure}
   
       \caption{}
     \label{fig:weakly-closed-8}
 \end{figure}

Consider the following ordering on the vertices
$$x_{n'}<_{lex}x_{(n-1)'}<_{lex} \cdots<_{lex}x_{1'}<_{lex}x_{0'}<_{lex}x_n<_{lex}\cdots <_{lex} x_1.$$
Let $<_{lex}$ denote the graded lexicographic ordering induced by this total ordering on the vertices. 
For each $1\leq i\leq n$, let $f_i$ denote the inner minor whose diagonal or anti-diagonal vertices are $\{i,i'\}$. Note that with respect to the graded lexicographic ordering, $\ini_{<_{lex}} (f_i)=x_ix_{i'}$ (in all the figures the initial terms are also drawn by diagonal or anti-diagonal lines). Also, by the construction of the total ordering, it is evident that $\gcd(\ini_{<_{lex}}(f_i), \ini_{<_{lex}}(f_j))=1$ for all $i\neq j$. Hence $f_1,\ldots,f_n$ forms a regular sequence. Since $\height(I_\mathcal{P})=\rank(\mathcal{P})=n$ the proof is complete. 
\end{proof}
\begin{figure}[h]%
		\centering
  		\begin{tikzpicture}[scale=.8]
			\draw[fill opacity=.5] (1,0) rectangle (5,1);
			\draw[fill opacity=.5] (1,1) rectangle (2,2);
			\draw[fill opacity=.5] (4,1) rectangle (5,2);
			\draw[fill opacity=.5] (0,1) rectangle (1,5);
			\draw[fill opacity=.5] (1,5) rectangle (3,6);
			\draw[fill opacity=.5] (2,6) rectangle (6,7);
			\draw[fill opacity=.5] (6,2) rectangle (7,5);
			\draw[fill opacity=.5] (5,4) rectangle (6,6);
			\draw[fill opacity=.5] (5,2) rectangle (6,3);
			\draw[fill opacity=.5] (1,4) rectangle (2,5);
			\draw (2,0)--(2,1);
			\draw (3,0)--(3,1);
			\draw (4,0)--(4,1);
			\draw (5,1)--(5,2);
			\draw (6,3)--(7,3);
			\draw (6,4)--(7,4);
			\draw (5,5)--(6,5);
			\draw (2,5)--(2,6);
			\draw (3,6)--(3,7);
			\draw (4,6)--(4,7);
			\draw (5,6)--(5,7);
			\draw (0,2)--(1,2);
			\draw (0,3)--(1,3);
			\draw (0,4)--(1,4);
			
			\draw (5,3)--(6,2);
			\node at (4.85,3.15)[ scale = 0.025cm] {$1$};
			\node at (6,1.75)[ scale = 0.025cm] {$1'$};
			
			\draw (6,3)--(7,2);
			\node at (5.85,3.25)[ scale = 0.025cm] {$2'$};
			\node at (7,1.75)[ scale = 0.025cm] {$2$};
			
			\draw (6,4)--(7,3);
			\node at (5.75,3.75)[ scale = 0.025cm] {$3'$};
			\node at (7.25,3)[ scale = 0.025cm] {$3$};
			
			\draw (6,5)--(7,4);
			\node at (6.25,5.25)[ scale = 0.025cm] {$5'$};
			\node at (7.25,4)[ scale = 0.025cm] {$5$};
			
			\draw (5,4)--(7,5);
			\node at (4.85,3.85)[ scale = 0.025cm] {$4$};
			\node at (7.25,5.25)[ scale = 0.025cm] {$4'$};
			
			\draw (5,5)--(6,6);
			\node at (4.85,4.85)[ scale = 0.025cm] {$6$};
			\node at (6.25,6)[ scale = 0.025cm] {$6'$};

			\draw (5,6)--(6,7);
			\node at (4.75,5.75)[ scale = 0.025cm] {$7'$};
			\node at (6.25,7)[ scale = 0.025cm] {$7$};
			
			\draw (4,6)--(5,7);
			\node at (3.85,5.75)[ scale = 0.025cm] {$8'$};
			\node at (5,7.25)[ scale = 0.025cm] {$8$};
			
			\draw (3,6)--(4,7);
			\node at (3.25,5.75)[ scale = 0.025cm] {$9'$};
			\node at (4,7.25)[ scale = 0.025cm] {$9$};
			
			\draw (2,6)--(3,7);
			\node at (1.75,6.25)[ scale = 0.025cm] {$11'$};
			\node at (3,7.25)[ scale = 0.025cm] {$11$};
			
			\draw (2,7)--(3,5);
			\node at (3.25,5)[ scale = 0.025cm] {$10$};
			\node at (1.75,7.25)[ scale = 0.025cm] {$10'$};
			
			\draw (1,6)--(2,5);
			\node at (2.35,4.75)[ scale = 0.025cm] {$12'$};
			\node at (0.75,6.25)[ scale = 0.025cm] {$12$};
			
			\draw (1,5)--(2,4);
			\node at (2.25,3.75)[ scale = 0.025cm] {$13$};
			\node at (0.65,5.25)[ scale = 0.025cm] {$13'$};
			
			\draw (0,5)--(1,4);
			\node at (-0.25,5.25)[ scale = 0.025cm] {$14$};
			\node at (1.35,3.75)[ scale = 0.025cm] {$14'$};
			
			\draw (0,4)--(1,3);
			\node at (-0.25,4.25)[ scale = 0.025cm] {$15$};
			\node at (1.35,3)[ scale = 0.025cm] {$15'$};
			
			\draw (0,3)--(1,2);
			\node at (-0.25,3.25)[ scale = 0.025cm] {$16$};
			\node at (1.35,2.25)[ scale = 0.025cm] {$16'$};
			
			\draw (0,2)--(1,1);
			\node at (-0.25,2)[ scale = 0.025cm] {$18$};
			\node at (0.65,0.75)[ scale = 0.025cm] {$18'$};
			
			\draw (0,1)--(2,2);
			\node at (-0.25,1)[ scale = 0.025cm] {$17$};
			\node at (2.35,2.15)[ scale = 0.025cm] {$17'$};
			
			\draw (2,1)--(1,0);
			\node at (0.75,-0.15)[ scale = 0.025cm] {$19$};
			\node at (2.35,1.25)[ scale = 0.025cm] {$19'$};
			
			\draw (3,1)--(2,0);
			\node at (1.75,-0.25)[ scale = 0.025cm] {$20$};
			\node at (3,1.25)[ scale = 0.025cm] {$20'$};
			
			\draw (4,1)--(3,0);
			\node at (2.75,-0.25)[ scale = 0.025cm] {$21$};
			\node at (3.75,1.25)[ scale = 0.025cm] {$21'$};
			
			\draw (5,1)--(4,0);
			\node at (3.75,-0.25)[ scale = 0.025cm] {$23$};
			\node at (5.35,1.25)[ scale = 0.025cm] {$23'$};

			\draw (5,0)--(4,2);
			\node at (4,2.25)[ scale = 0.025cm] {$22'$};
			\node at (5.35,-0.25)[ scale = 0.025cm] {$22$};
			
			\node at (4.75,2.25)[ scale = 0.025cm] {$0'$};
		\end{tikzpicture}
		\caption{}%
	\label{example-weakly-closed}
\end{figure}

\begin{example}
 We consider a weakly closed path polyomino with labeling of vertices as in Figure~\ref{example-weakly-closed}. 
	Consider the graded lexicographic order induced by
 $$x_{23'}<_{lex}\cdots <_{lex}x_{1'}<_{lex}x_{0'}<_{lex}x_{23}<_{lex} \cdots<_{lex}x_1.$$
Then one can see that  with respect to this ordering, the polyomino ideal is of K\"onig type, where  diagonal or anti-diagonal lines are drawn to denote $\ini_{<_{lex}}(f_i)$  for all $i=1,\cdots, 23$.
	\end{example}

 \begin{corollary}\label{thm:weaklyclosedknutson}
Let $\mathcal{P}: A_0,\cdots,A_{n-1}$ be a weakly closed path polyomino.
 Then $I_{\mathcal{P}}$ is a Knutson ideal.	    
\end{corollary}
\begin{proof}
In the proof of the above theorem, we observed that for each cell $A_i$, we obtain a generator $f_i$ of $I_\MP$ corresponding to an inner interval of $\MP$ such that $\gcd(\ini_{<_{lex}}(f_i), \ini_{<_{lex}} (f_j))=1$ for all $i\neq j$. 
    Consider $f=\prod f_i$, then $\ini_{<_{lex}}(f)$ is square-free.     
    We show that $I_{\mathcal{P}}\in \mathcal{C}_f$. 
    Let $\mathcal{P}_1=\mathcal{P}\setminus A_0$ and $\mathcal{P}_2=\mathcal{P}\setminus A_{n-1}$. 
    Since $\MP_1$ and $\MP_2$ are simple polyominoes, $I_{\MP_1}$ and $I_{\MP_2}$ are prime ideals.
    Observe that $I_{\MP_1}$ (resp. $I_{\MP_2}$) is a minimal prime ideal of $(f_1,f_2,\ldots, f_{n-1})$ (resp. $(f_0,f_2,\ldots, f_{n-2})$).  
    Since $(f_1,f_2,\ldots, f_{n-1}), (f_0,f_2,\ldots, f_{n-2}) \in \mathcal{C}_f$, we get that $I_{\MP_1}, I_{\MP_2} \in \mathcal{C}_f$.
    Also, note that $I_{\MP}=I_{\MP_1}+I_{\MP_2}$. Hence the proof follows.
\end{proof}	
\begin{Remark}
One can note that using similar total order on the vertices as of weakly closed path polyomino any thin path polyomino is K\"onig and hence the corresponding polyomino ideal is Knutson.
\end{Remark}

\section{Ladder polyominoes are Knutson}\label{sec:ladder}

In this section, we study a class of polyominoes that generalizes parallelogram polyominoes. 
We show that the polyomino ideals associated with these polyominoes are Knutson. 

\begin{figure}[h]{}
		\centering
		\begin{tikzpicture}[scale=.6]
			\draw[] (1,1)--(13,1)--(13,2)--(14,2)--(14,3)--(9,3)--(9,4)--(7,4)--(7,5)--(11,5)--(11,6)--(8,6)--(8,7)--(10,7)--(10,8)--(6,8)--(6,5)--(4,5)--(4,3)--(3,3)--(3,2)--(1,2)--(1,1);
			\draw[] (2,1)--(10,1);
			\draw[] (3,2)--(13,2);
			\draw[] (3,3)--(10,3);
			\draw[] (4,4)--(10,4);
			\draw[] (5,5)--(11,5);
			\draw[] (6,6)--(8,6);
			\draw[] (6,7)--(8,7);
			\draw[] (2,1)--(2,2);
			\draw[] (3,1)--(3,2);
			\draw[] (4,1)--(4,5);
			\draw[] (5,1)--(5,5);
			\draw[] (6,1)--(6,8);
			\draw[] (7,1)--(7,8);
			\draw[] (8,1)--(8,8);
			\draw[] (9,1)--(9,4);
			\draw[] (9,5)--(9,8);
			\draw[] (10,5)--(10,6);
			\draw[] (9,7)--(9,8);
			\draw[] (10,1)--(10,4);
			\draw[] (10,7)--(10,8);
			\draw[] (11,1)--(11,3);
			\draw[] (11,5)--(11,6);
			\draw[] (12,1)--(12,3);
			\draw[] (13,1)--(13,3);

          \filldraw[] (1,1) circle (0.5pt) node[anchor=east]  {$a_1$};
          \filldraw[] (1,2) circle (0.5pt) node[anchor=east]  {$a_2$};
          \filldraw[] (3,3) circle (0.5pt) node[anchor=east]  {$a_3$};
          \filldraw[] (4,4) circle (0.5pt) node[anchor=east]  {$a_4$};
          \filldraw[] (4,5) circle (0.5pt) node[anchor=east]  {$a_5$};
          \filldraw[] (6,6) circle (0.5pt) node[anchor=east]  {$a_6$};
          \filldraw[] (6,7) circle (0.5pt) node[anchor=east]  {$a_7$};
          \filldraw[] (6,8) circle (0.5pt) node[anchor=east]  {$a_8$};
	\end{tikzpicture}
        \caption{A ladder polyomino}%
	\label{fig:one-sided}
\end{figure}
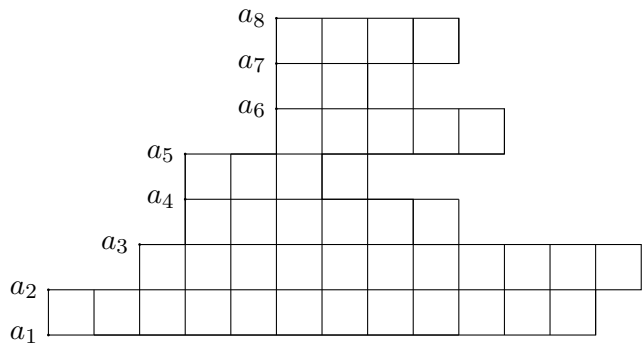

Let $\MP$ be a horizontally convex polyomino. 
A set of vertices $\{a_1 = (a_{1,1}, a_{1,2}), \ldots, a_r = (a_{r,1}, a_{r,2})\} \subset V(\MP)$, with $a_{1,2}<a_{2,2}<\ldots< a_{r,2}$ is said to be the set of \emph{left-most vertices} of $\MP$, if for any vertex $v = (v_1, v_2) \in V(\MP)$, either $v = a_i$ for some $i$, or there exists a $j,$ $1\leq j \leq r$ such that $v_2 = a_{j,2}$, and $a_{j,1} < v_1$.
A horizontally convex polyomino $\MP$ with the set of left-most vertices $\{a_1, a_2, \ldots, a_r\}$ is called a \emph{ladder polyomino} if for each $2\leq i \leq r,$ $a_i$ is an upper-left vertex of a cell in $\MP.$
An example of ladder polyomino is shown in Figure~\ref{fig:one-sided} with $\{a_1,\ldots,a_8\}$ as the set of left-most vertices.
Note that in a ladder polyomino, $a_{1,1} \leq a_{2,1} \leq \ldots \leq a_{r,1}$ and $a_{i,2}=a_{i-1,2}+1$ for all $2\leq i\leq r$.
Also, note that a parallelogram polyomino is a ladder polyomino.

\begin{theorem}\label{thm:ladderpolyomino}
Let $\MP$ be a ladder polyomino. Then, $I_{\MP}$ is a Knutson ideal.
\end{theorem}

Before we proceed with the proof of Theorem~\ref{thm:ladderpolyomino}, we will establish several lemmas that are essential to the proof. We begin by introducing some notations that will be utilized not only in the proof of this theorem but also throughout the remainder of the article.

\begin{discussionbox}
    \label{setup:ladder}
Let $\MP$ be a collection of cells. 
Partition the vertices of $\MP$ into sets $V_1, \ldots, V_s$ which satisfy following conditions: for each $1\leq k\leq s$, let $V_k := \{v_{k,1},\ldots,v_{k,n_{k}}\}$, then
     \begin{asparaenum}
      \item
    for all  $1\leq i < n_k$, the vertices $v_{k,i}$ and $v_{k,{i+1}}$ are upper-left and lower-right vertices of a cell in $\MP$ respectively.
       \item
       $v_{k,1}$ (resp. $v_{k,n_{k}}$ ) is not a lower-right (resp. upper-left) vertex of any cell in $\MP$.
       \item 
       for $k<l,$ write $v_{k,1}=(v_1,v_2)\in V_k$ and $v_{l,1}=(w_1,w_2)\in V_l$, then $v_1+v_2<w_1+w_2$,  or $v_1+v_2 =w_1+w_2$ and $v_1<w_1$.
      \end{asparaenum}

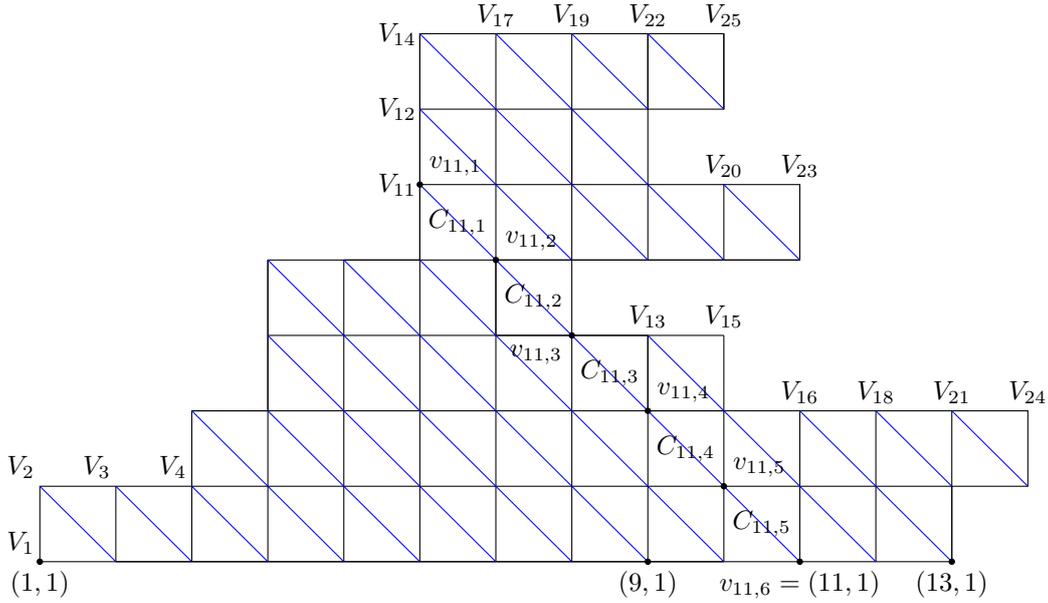
\begin{figure}[h]{}
\small
		\centering
		\begin{tikzpicture}[scale=1]
			\draw[] (1,1)--(13,1)--(13,2)--(14,2)--(14,3)--(9,3)--(9,4)--(7,4)--(7,5)--(11,5)--(11,6)--(8,6)--(8,7)--(10,7)--(10,8)--(6,8)--(6,5)--(4,5)--(4,3)--(3,3)--(3,2)--(1,2)--(1,1);
			\draw[] (2,1)--(10,1);
			\draw[] (3,2)--(13,2);
			\draw[] (3,3)--(10,3);
			\draw[] (4,4)--(10,4);
			\draw[] (5,5)--(11,5);
			\draw[] (6,6)--(8,6);
			\draw[] (6,7)--(8,7);
			\draw[] (2,1)--(2,2);
			\draw[] (3,1)--(3,2);
			\draw[] (4,1)--(4,5);
			\draw[] (5,1)--(5,5);
			\draw[] (6,1)--(6,8);
			\draw[] (7,1)--(7,8);
			\draw[] (8,1)--(8,8);
			\draw[] (9,1)--(9,4);
			\draw[] (9,5)--(9,8);
			\draw[] (10,5)--(10,6);
			\draw[] (9,7)--(9,8);
			\draw[] (10,1)--(10,4);
			\draw[] (10,7)--(10,8);
			\draw[] (11,1)--(11,3);
			\draw[] (11,5)--(11,6);
			\draw[] (12,1)--(12,3);
			\draw[] (13,1)--(13,3);
			\draw[color=blue] (1,2)--(2,1);
			\draw[color=blue] (2,2)--(3,1);
			\draw[color=blue] (3,2)--(4,1);
			\draw[color=blue] (3,3)--(5,1);
			\draw[color=blue] (4,3)--(6,1);
			\draw[color=blue] (4,4)--(7,1);
			\draw[color=blue] (4,5)--(8,1);
			\draw[color=blue] (5,5)--(9,1);
			\draw[color=blue] (6,5)--(10,1);
			\draw[color=blue] (6,6)--(11,1);
			\draw[color=blue] (6,7)--(8,5);
			\draw[color=blue] (9,4)--(12,1);
			\draw[color=blue] (6,8)--(9,5);
			\draw[color=blue] (11,3)--(13,1);
			\draw[color=blue] (7,8)--(10,5);
			\draw[color=blue] (13,3)--(14,2);
			\draw[color=blue] (12,3)--(13,2);
			\draw[color=blue] (8,8)--(9,7);
			\draw[color=blue] (9,8)--(10,7);
			\draw[color=blue] (10,6)--(11,5);
			
                \filldraw[] (1,1) circle (1pt) node[anchor=north]  {$(1,1)$};
                \filldraw[] (9,1) circle (1pt) node[anchor=north]  {$(9,1)$};
                \filldraw[] (13,1) circle (1pt) node[anchor=north]  {$(13,1)$};
			\node[] at (.75,1.25) () {$V_{1}$};
			\node[] at (.75,2.25) () {$V_{2}$};
			\node[] at (1.75,2.25) () {$V_{3}$};
			\node[] at (2.75,2.25) () {$V_{4}$};
			\node[] at (10,8.25) () {$V_{25}$};
			\node[] at (9,8.25) () {$V_{22}$};
			\node[] at (11,6.25) () {$V_{23}$};
			\node[] at (10,6.25) () {$V_{20}$};
			\node[] at (14,3.25) () {$V_{24}$};
			\node[] at (5.7,6) () {$V_{11}$};
			\node[] at (5.7,7) () {$V_{12}$};
			\node[] at (5.7,8) () {$V_{14}$};
			\node[] at (7,8.25) () {$V_{17}$};
			\node[] at (8,8.25) () {$V_{19}$};
			\node[] at (13,3.25) () {$V_{21}$};
			\node[] at (12,3.25) () {$V_{18}$};
			\node[] at (11,3.25) () {$V_{16}$};
			\node[] at (10,4.25) () {$V_{15}$};
			\node[] at (9,4.25) () {$V_{13}$};

        \filldraw[] (6,6) circle (1pt) node[anchor=south west]  {$v_{11,1}$};
        \filldraw[] (7,5) circle (1pt) node[anchor=south west]  {$v_{11,2}$};
        \filldraw[] (8,4) circle (1pt) node[anchor=north east]  {$v_{11,3}$};
        \filldraw[] (9,3) circle (1pt) node[anchor=south west]  {$v_{11,4}$};
        \filldraw[] (10,2) circle (1pt) node[anchor=south west]  {$v_{11,5}$};
        \filldraw[] (11,1) circle (1pt) node[anchor=north]  {$v_{11,6}=(11,1)$};
        
        \filldraw[] (6.5,5.5) circle (0pt) node[anchor=center]  {$C_{11,1}$};
        \filldraw[] (7.5,4.5) circle (0pt) node[anchor=center]  {$C_{11,2}$};
        \filldraw[] (8.5,3.5) circle (0pt) node[anchor=center]  {$C_{11,3}$};
        \filldraw[] (9.5,2.5) circle (0pt) node[anchor=center]  {$C_{11,4}$};
        \filldraw[] (10.5,1.5) circle (0pt) node[anchor=center]  {$C_{11,5}$};
	\end{tikzpicture}
 \caption{A ladder polyomino}%
	\label{fig:labeledone-sided}
\end{figure}

In loose words, $V_i$'s are the sets of vertices on anti-diagonal lines. For $i < j$, either $V_j$ is on the anti-diagonal line to the right of $V_i$, or if they reside on the same anti-diagonal line, then $V_j$ is positioned lower than $V_i$.
An example is given in Figure~\ref{fig:labeledone-sided} emphasizing the set $V_{11}$.

To each $V_k$, we now associate a set of cells $\MC_k$ in the following way: If  $V_k$ is a singleton, then we define $\mathcal{C}_k$ to be the empty-set. Otherwise, define $\mathcal{C}_k := \{C_{k,1}, C_{k,2}, \ldots, C_{k,{n_k - 1}}\}$ as the set of cells in $\MP$ where, for each $1\leq i<n_k-1$, $v_{k,i}$ and $v_{k,{i+1}}$ are the anti-diagonal vertices of the cell $C_{k,i}$. The set $\MC_{11}$ is labeled in Figure~\ref{fig:labeledone-sided}. Let $\MP_{k}$ denote the sub-collection of cells of $\MP$, defined as:
	\[
	\MP_{k} := \bigcup_{i \leq k} \mathcal{C}_i.
	\]

For any $1\leq k\leq s$, suppose that $v_{k,1} = (i_k+1, j_k + n_k)$ for some $i_k,j_k\in \naturals$, 
then $v_{k,l} = (i_k+l, j_k + n_k+1-l)$ for all $1\leq l\leq n_k$.
To each $V_k$, we associate a matrix $M_k=[m_{ab}]_{1\leq a,b\leq n_k}$ of order $n_k$ with the following entries: for $1\leq a,b\leq n_k$,
\[ m_{ab} = \begin{cases} 
          x_{i_k+a, j_k + b} & \text{if}\ (i_k+a, j_k + b)\in V(\MP), \\
         0 & \text{else.} 
       \end{cases}
    \]
Note that the diagonal entries of the matrix $M_k$ are the variables corresponding to the vertices in $V_k$. 
The following matrix represents $M_{13}$ associated with $V_{13}$ as in Figure~\ref{fig:labeledone-sided}.

\[
\begin{bmatrix}
    x_{9,4}       & x_{10,4} & 0 & 0 \\
    x_{9,3}       & x_{10,3} & x_{11,3} & x_{12,3} \\
    x_{9,2}       & x_{10,2} & x_{11,2} & x_{12,2} \\
    x_{9,1}       & x_{10,1} & x_{11,1} & x_{12,1} \\
\end{bmatrix}.
\]
Let $f_{k}$ denote the determinant of the matrix $M_k$. 
Define 
\[f := \prod_{k=1}^{s} f_{k}.
\]
 We define a total order $<$ on the variables of $S$ as follows: $x_{i,j} < x_{k,l}$ if $i < k$ or $i = k$ and $j < l$. We apply the graded reverse lexicographic order induced by the order $<$. Note that the initial term of $x_{a} x_{b} - x_{c} x_{d},$ associated with the inner interval $[a,b]$ of $\MP$ with anti-diagonal vertices $c$ and $d$, is $x_cx_d$ with respect to the order defined above.
\end{discussionbox}

\begin{lemma}
	\label{lem:one_sided_initial}
	Under the notations of Discussion~\ref{setup:ladder},
	\[
	\ini_{<}(f) = \prod_{v \in V(\MP)} x_v.
	\]
\end{lemma}

\begin{proof}
    It suffices to show that $\ini_{<}(f_{k}) = \prod_{v \in V_k} x_v$ for every $1\leq k\leq s$.
    If $|V_k|= 1$, then we are done.
    So we may assume that $|V_k|\geq 2$.
    Recall that $v_{k,l} = (i_k+l, j_k + n_k+1-l)$ for all $1\leq l\leq n_k$ and for some $i_k,j_k\in \naturals$.
    The matrix $M_k$ has the following form 
\[
    \begin{bmatrix}
    x_{i_k+1, j_k + n_k}     & * & \dots & * \\
    *    & x_{i_k+2, j_k + n_k-1} & \dots & * \\
    \vdots      & \vdots & \ddots & \vdots \\
    *       & * & \dots & x_{i_k + n_k, j_k+1} \\
\end{bmatrix}.
\]
Therefore, 
\[
f_k = \sum_{\sigma \in S_{n_k}} (-1)^{\sgn(\sigma)} x_{i_k+1, j_k + \sigma(1)} x_{i_k+2, j_k + \sigma(2)} \cdots x_{i_k + n_k, j_k + \sigma(n_k)},
\]
where $S_{n_k}$ denotes the symmetric group on $n_k$ symbols, and $\sgn(\sigma)$ denotes the sign of the permutation $\sigma$, which is $1$ if $\sigma$ is even and $-1$ if $\sigma$ is odd.
Note that some of the monomials in the summation might not be present.

It is enough to show that
	\[
	x_{i_k+1, j_k + n_k} x_{i_k+2, j_k + n_k - 1} \cdots x_{i_k + n_k, j_k + 1} > x_{i_k+1, j_k + \sigma(1)} x_{i_k+2, j_k + \sigma(2)} \cdots x_{i_k + n_k, j_k + \sigma(n_k)}
	\]
	whenever the right-hand side is non-zero and $\sigma$ is not the permutation
	\[
	\begin{pmatrix}
		1 & 2 & \cdots & n_k \\
		n_k & n_k - 1 & \cdots & 1
	\end{pmatrix}.
	\]
 This follows immediately since we are considering the graded reverse lexicographic order. 
\end{proof}

Before proving Theorem~\ref{thm:ladderpolyomino}, we need a technical result (Lemma~\ref{lem:detf_k}). To prove Lemma~\ref{lem:detf_k} we develop a combinatorial method which we discuss below.

\begin{setup}
	\label{setup_chi}
	Let $\chi(n, l, \sigma)$ be the following property:
	for an $n \geq 2$, an integer $l$ with $2\leq l \leq n$, and a permutation $\sigma \in S_n$ (the symmetric group on $n$ symbols), either one of the following holds:
	\begin{enumerate}[label=\textbf{\ref{setup_chi}.\arabic*}]
		\item \label{chi_a}
        There exists an integer $i$ with $i < l$, such that $l + \max\{\sigma(i), \sigma(l)\} \leq n + 2$.
		
		\item \label{chi_b}
        There exist two integers $i, j$ with $1 \leq i < j \leq n$ such that $j + \max\{\sigma(i), \sigma(j)\} \leq n + 1$.
	\end{enumerate}
\end{setup}

\begin{lemma}
		\label{lemma_chi_n}
		For every $n \geq 2$, $2\leq l \leq n$ and $\sigma\in S_n$, $\chi(n, l, \sigma)$ holds true.
\end{lemma}

\begin{proof}
    We proceed by induction on $n$. 
    If $n=2,$ then $l=2$ and \ref{chi_a} holds for every $\sigma\in S_2$.
    So we may assume that $n\geq 3$ and $\chi(n-1, l, \sigma)$ holds for every $2 \leq l \leq n-1$ and $\sigma\in S_{n-1}$.
    Fix a $\sigma\in S_n$ and $2\leq l\leq n$.
    Note that if $l = 2$, then $2 + \max\{\sigma(1), \sigma(2)\} \leq n + 2$; thus \ref{chi_a} holds.
    So we may assume that $l \geq 3.$
    Let $l' = l - 1$, and $n_0$ with $1\leq n_0\leq n$ be such that $\sigma(n_0)=n$.
    Define a permutation $\widetilde{\sigma} \in S_{n-1}$ by
    \[
	\widetilde{\sigma}(m) = 
		\begin{cases}
			\sigma(m) & \text{if } m < n_0 \\
			\sigma(m+1) & \text{if } m \geq n_0.
		\end{cases}
    \]
    
    Since $\chi(n-1, l', \widetilde{\sigma})$ holds by induction hypothesis, either \ref{chi_a} or \ref{chi_b} holds for $n-1, l'$ and $\widetilde{\sigma}$.    
     First, assume that \ref{chi_a} holds. 
     Then, there exists an $i$ with $i < l'$, such that $l' + \max\{\widetilde{\sigma}(i), \widetilde{\sigma}(l')\} \leq n + 1$. 
The following cases are possible:
\begin{asparaenum}
       \item
           When $i<l'<n_0$. Then, $\widetilde{\sigma}(i)= \sigma(i), \widetilde{\sigma}(l')=\sigma(l')$. 
           So, $l' + \max\{\sigma(i), \sigma(l')\} \leq n+1$; thus \ref{chi_b} holds for $\sigma$ and $l$ with respect to $i$ and $l'$.
       
       \item
           When $i<n_0\leq l'$. Then, $\widetilde{\sigma}(i)=\sigma(i), \widetilde{\sigma}(l')= \sigma(l'+1)$.
           So, $l + \max\{\sigma(i), \sigma(l)\} \leq n+2$; thus \ref{chi_a} holds for $\sigma$ and $l$.

        \item
            When $n_0\leq i< l'$. Then, $\widetilde{\sigma}(i)= \sigma(i+1), \widetilde{\sigma}(l')=\sigma(l'+1)$.
           So, $l + \max\{\sigma(i+1), \sigma(l)\} \leq n+2$; thus \ref{chi_a} holds for $\sigma$ and $l$ with respect to $i+1$.
       \end{asparaenum}
		
    On the other hand, suppose that \ref{chi_b} holds for $n-1, l'$ and $\widetilde{\sigma}$ . Then, there exist $1 \leq i < j \leq n-1$ such that $j + \max\{\widetilde{\sigma}(i), \widetilde{\sigma}(j)\} \leq n$.
    The following cases are possible:
\begin{asparaenum}
       \item
           When $i<j<n_0$. Then, $\widetilde{\sigma}(i)=\sigma(j), \widetilde{\sigma}(j)= \sigma(j)$. 
           So, $j + \max\{\sigma(i), \sigma(j)\} \leq n < n+1$; thus \ref{chi_b} holds for $\sigma$ and $l$.
       
       \item
           When $i<n_0\leq j$. Then, $\widetilde{\sigma}(i)=\sigma(i), \widetilde{\sigma}(j)= \sigma(j+1)$.
           So, $(j+1) + \max\{\sigma(i), \sigma(j+1)\} \leq n+1$; thus \ref{chi_b} holds for $\sigma$ and $l$ with respect to $i$ and $j+1$.

        \item
            When $n_0\leq i< j$. Then, $\widetilde{\sigma}(i)=\sigma(i+1), \widetilde{\sigma}(j)= \sigma(j+1)$.
           So, $(j+1) + \max\{\sigma(i+1), \sigma(j+1)\} \leq n+1$; thus \ref{chi_b} holds for $\sigma$ and $l$ with respect to $i+1$ and $j+1$.
       \end{asparaenum}
\end{proof}

For some fixed $n\geq 2, 2 \leq l \leq n$ and  $\sigma \in S_n$, let 
\[ S(n,l,\sigma) := \{(i,j) : i < j,\ \text{either}\ j=l\ \text{and \ref{chi_a} holds for}\ i,j \ \text{or}\ \text{\ref{chi_b} holds for}\ i,j\}.
 \]
 By Lemma \ref{lemma_chi_n}, $S(n,l,\sigma)$ is non-empty.
Let $(i_\sigma,j_\sigma)$ be the minimum element of $S(n,l,\sigma)$ under the graded lexicographic order on $\naturals^2$.

\begin{Lemma}\label{lem:minimal}
    Under the notations as above, $(i_\sigma, j_\sigma)$ is the minimum element of $S(n,l, \sigma\circ (i_\sigma j_\sigma)),$
    where $(i_\sigma j_\sigma)$ denotes the transposition of $S_n$.
\end{Lemma}

\begin{proof}
    Let $\sigma' = \sigma\circ (i_\sigma j_\sigma).$
    Since $\sigma'(i_\sigma) = \sigma(j_\sigma)$ and $\sigma'(j_\sigma) = \sigma(i_\sigma)$, we get that $(i_\sigma,j_\sigma)\in S(n,l,\sigma').$ 
    Let $(i,j)$ be a pair with $(i,j)<(i_\sigma, j_\sigma)$ under the graded lexicographic order on $\naturals^2$. 
    The following cases are possible: either $\{i,j\}\cap \{i_\sigma, j_\sigma\}=\varnothing$ or $\#(\{i,j\}\cap \{i_\sigma, j_\sigma\})=1$. As $(i,j)<(i_\sigma, j_\sigma),$ in the later case, either $i<i_\sigma$ and $j=j_\sigma$, or $i=i_\sigma$ and $j<j_\sigma$. We proceed in all cases.

    \begin{asparaenum}
       \item
           When $\{i,j\}\cap \{i_\sigma, j_\sigma\}=\varnothing$. 
           If $(i,j)\in S(n,l,\sigma'),$ then the property $\chi(n, l, \sigma')$ holds for the pair $(i,j)$.
           Since $\sigma'(i) = \sigma(i)$ and $\sigma'(j) = \sigma(j)$, we get that $\chi(n, l, \sigma)$ holds for the pair $(i,j)$.
           Therefore, this case is not possible.
       
       \item 
           When $i<i_\sigma$ and $j=j_\sigma$. Since $(i_\sigma, j_\sigma)\in S(n,l,\sigma)$, we have $j_\sigma+\sigma(i_\sigma) \leq n + 2$ if $j_\sigma= l$, and $\leq n+1$ else; thus $i_\sigma+\sigma(i_\sigma) \leq n + 1$ because $i_\sigma<j_\sigma.$ 
           On the other hand, since $(i,i_\sigma)<(i_\sigma, j_\sigma)$, we have $(i,i_\sigma)\notin S(n,l,\sigma)$. 
           Thus, $i_\sigma + \max\{\sigma(i), \sigma(i_\sigma)\} \nleq n + 2$ if $i_\sigma= l$, and $\nleq n+1$ else. 
           Since $i_\sigma +\sigma(i_\sigma) \leq n + 1$, we get that  $i_\sigma + \sigma(i) \nleq n + 2$ if $i_\sigma= l$, and $\nleq n+1$ else.   
           Since, $\sigma(i)= \sigma'(i),$ we get $(i,j)\notin S(n,l,\sigma')$.

        \item
            When $i=i_\sigma$ and $j<j_\sigma$. Since $(i_\sigma, j_\sigma)\in S(n,l,\sigma)$, we have $j_\sigma+\sigma(i_\sigma) \leq n + 2$; thus, $j +\sigma(i_\sigma) \leq n + 1$ because $j<j_\sigma$.
            On the other hand, since $(i,j)<(i_\sigma, j_\sigma)$, we have $(i,j)\notin S(n,l,\sigma)$. 
            Thus, $j + \max\{\sigma(i), \sigma(j)\} \nleq n + 2$ if $j= l$, and $\nleq n+1$ else.
            Since $j +\sigma(i=i_\sigma) \leq n + 1$, we get that  $j + \sigma(j) \nleq n + 2$ if $j= l$, and $\nleq n+1$ else.
            Since $\sigma(j)= \sigma'(j),$ we get $(i,j)\notin S(n,l,\sigma')$.
       \end{asparaenum}
       Hence the proof.
\end{proof}

Following the above notations for an $n\geq2 $ and $2\leq l\leq n$, 
let $\sigma' = \sigma\circ (i_\sigma j_\sigma)$, then by Lemma~\ref{lem:minimal}, $(i_{\sigma'}, j_{\sigma'}) =(i_\sigma, j_\sigma)$.  Thus if ${\sigma}_1\circ (i_{{\sigma}_1} j_{{\sigma}_1}) = {\sigma}_2\circ (i_{{\sigma}_2} j_{{\sigma}_2})=\sigma'$ (say), for some ${\sigma}_1,{\sigma}_2\in S_n$,  we have  $(i_{{\sigma}_1}, j_{{\sigma}_1})= (i_{\sigma'}, j_{\sigma'}) =(i_{{\sigma}_2}, j_{{\sigma}_2})$. Hence, ${\sigma}_1 = {\sigma}_2$. Note that if $\sigma$ varies over all even permutations, then by what we discussed above, the set $\{\sigma\circ (i_\sigma j_\sigma): \sigma\in S_n \textrm{ a even permutation }\}$ is the set of all odd permutations. 
Thus, we can partition the set $S_n =  \bigsqcup_{\sigma \in \Lambda} \{\sigma, \sigma\circ(i_\sigma j_\sigma)\},$ where $\Lambda$ is the set of all even permutation of $S_n$.

  \begin{figure}[h]
		\centering
		\begin{tikzpicture}[scale=1]
		\draw[thick] 
        (1,1)--(4,1) (5,1)--(7,1) 
        (1,2)--(4,2) (5,2)--(7,2) 
        (1,3)--(4,3) (5,3)--(6,3)

        (1,4)--(4,4)
        (1,5)--(4,5)
        (1,6)--(3,6)
        (1,7)--(2,7)
       
        (1,1)--(1,3) (1,4)--(1,7)
        (2,1)--(2,3) (2,4)--(2,7)
        (3,1)--(3,3) (3,4)--(3,6)
        (4,1)--(4,3) (4,4)--(4,5)
        
        (5,1)--(5,3)
        (6,1)--(6,3)
        (7,1)--(7,2)

        ;
        \draw[dotted, thick]
        (4,1)--(5,1)
        (4,2)--(5,2)
        (4,3)--(5,3)
        (4,4)--(5,4)

        (1,4)--(1,3)
        (2,4)--(2,3)
        (3,4)--(3,3)
        (4,4)--(4,3)
        (5,3)--(5,4)

        (7,1)--(8,1)--(8,2)--(7,2)
        (8,2)--(9,2)--(9,3)--(6,3)
        
        (2,7)--(4,7)--(4,6)--(3,6)
        (4,6)--(6,6)--(6,5)--(4,5)
        (5.5,5)--(5.5,4)--(5,4)
        (5.5,4)--(7,4)--(7,3)
        ;
        \fill[fill=gray, fill opacity=0.3] (1,1)--(7,1)--(7,2)--(6,2)--(6,3)--(5,3)--(5,4)--(4,4)--(4,5)--(3,5)--(3,6)--(2,6)--(2,7)--(1,7)--(1,1);

        \filldraw[] (1.5,6.5) circle (0pt) node[anchor=center]  {$C_{k,1}$};
        \filldraw[] (2.5,5.5) circle (0pt) node[anchor=center]  {$C_{k,2}$};
        \filldraw[] (3.5,4.5) circle (0pt) node[anchor=center]  {$C_{k,3}$};
        \filldraw[] (5.5,2.5) circle (0pt) node[anchor=center]  {$C_{k,k-2}$};
        \filldraw[] (6.5,1.5) circle (0pt) node[anchor=center]  {$C_{k,k-1}$};

        \filldraw[] (1.5,5.5) circle (0pt) node[anchor=center]  {$C_{k-1,1}$};
        \filldraw[] (2.5,4.5) circle (0pt) node[anchor=center]  {$C_{k-1,2}$};
        \filldraw[] (1.5,1.5) circle (0pt) node[anchor=center]  {$C_{2,1}$};

        \filldraw[] (1,1) circle (0.5pt) node[anchor=north]  {$(1,1)$};
        \filldraw[] (7,1) circle (0.5pt) node[anchor=north]  {$(k,1)$};
        \filldraw[] (1,7) circle (0.5pt) node[anchor=east]  {$(1,k)$};
         \filldraw[] (1,2) circle (0.5pt) node[anchor=east]  {$(1,2)$};
          \filldraw[] (1,6) circle (0.5pt) node[anchor=east]  {$(1,k-1)$};

		\end{tikzpicture}
		\caption{}%
	\label{fig:lemma}
\end{figure}

\begin{Lemma}\label{lem:detf_k}
    Let $\MP$ be a horizontally convex polyomino as shown in Figure~\ref{fig:lemma} with coloured sub-polyomino.
    Under the notations of Discussion~\ref{setup:ladder}, if $\MC_k\neq \varnothing,$ then $f_k\notin I_{\MP_{k-1}}$ and $f_k\in I_{\MP_{k-1}\cup \{C\}}$ for all $C\in \MC_k.$
\end{Lemma}

\begin{proof}
Assume the notations of Discussion~\ref{setup:ladder}. Note that $V_{1} = (1,1), V_{2} = \{(1,2),(2,1)\},$
$V_{k-1}=\{(1,k-1),(2,k-2),\ldots,(k-1,1)\}$ and $V_k=\{(1,k),(2,k-1),\ldots,(k,1)\}.$

   Recall  from the proof of Lemma~\ref{lem:one_sided_initial},
    \[
f_k= \sum_{\sigma \in S_{k}} (-1)^{\sgn(\sigma)} x_{1, \sigma(1)} x_{2,\sigma(2)} \cdots x_{k,\sigma(k)},
\]
where $S_{n}$ denotes the symmetric group on $n$ symbols, and $\sgn(\sigma)$ denotes the sign of the permutation $\sigma$, which is $1$ if $\sigma$ is even and $-1$ if $\sigma$ is odd. 
Note that some of the monomials in the summation might not be present.

If $f_k\in I_{\MP_{k-1}}$, then $f_k = \sum_{i=1}^{r}h_ig_i,$ where $h_i\in S$ and $g_i\in I_{\MP_{k-1}}$ a generator corresponding to an inner interval of $\MP_{k-1}$ for all $1\leq i\leq r$.
So, every non-zero monomial of $f_k$ must contain at least two variables corresponding to diagonal or anti-diagonal vertices of an inner interval of $\MP_{k-1}$. 
Note that $x_{1, k} x_{2,k-1} \cdots x_{k,1}$ is a non-zero monomial in $f_k$, 
and it does not contain variables corresponding to diagonal or anti-diagonal vertices of any inner interval of $\MP_{k-1}$.
Thus, $f_k\notin I_{\MP_{k-1}}$. 

We now show that $f_k\in I_{\MP_{k-1}\cup \{C\}}$ for all $C\in \MC_k.$
Let $C_{k,l}\in \MC_k$. For $k$ and $l$, consider the partition of $S_k$ defined above.
Thus, we can rewrite $f_k$ as  
\begin{align*}
		f_k &= \sum_{\sigma \in \Lambda} \left(
		x_{1,\sigma(1)} x_{2, \sigma(2)} \cdots x_{k, \sigma(k)} 
		- x_{1,\sigma \circ (i_\sigma j_\sigma)(1)} x_{2,  \sigma \circ (i_\sigma j_\sigma)(2)} \cdots x_{k, \sigma \circ (i_\sigma j_\sigma)(k)} \right)\\
        &=\sum_{\sigma \in \Lambda} 
            \left(  
			x_{1,\sigma(1)} x_{2, \sigma(2)}\cdots \widehat{x_{{i_\sigma},\sigma(i_{\sigma})}} \cdots \widehat{x_{{j_\sigma},\sigma(j_{\sigma})}} \cdots x_{k, \sigma(k)} 
		\right)
		\left( 
		x_{{i_\sigma},\sigma(i_{\sigma})} x_{{j_\sigma},\sigma(j_{\sigma})} - x_{{i_\sigma},\sigma(j_{\sigma})} x_{{j_\sigma},\sigma(i_{\sigma})}
		\right),
	\end{align*}
 where $\widehat{x_{v}}$ indicates that the variable $x_v$ is not present in the product.

 Fix a $\sigma\in \Lambda$. Note that the pair $(i_{\sigma},j_{\sigma})\in S(k,l,\sigma)$,
 i.e., the pair $(i_{\sigma},j_{\sigma})$ satisfies \ref{chi_a} or \ref{chi_b}.
When $(i_{\sigma},j_{\sigma})$ satisfies~\ref{chi_b}, we get $i_{\sigma}< j_{\sigma}$ and $ j_{\sigma} + \max\{\sigma(i_{\sigma}), \sigma(j_{\sigma})\} \leq k + 1$. Consequently, it follows that
 \[
 (j_\sigma,\max\{\sigma(i_{\sigma}), \sigma(j_{\sigma})\})\in V(\MP_{k-1}).
 \]
 By the definition of $\MP_{k-1}$, we get that $(j_\sigma,\min\{\sigma(i_{\sigma}), \sigma(j_{\sigma})\})\in V(\MP_{k-1}).$
 Moreover, since $i_{\sigma}< j_{\sigma}$ and the way $\MP_{k-1}$ is defined, we get 
 \[
 (i_\sigma,\min\{\sigma(i_{\sigma}), \sigma(j_{\sigma})\}), (i_\sigma,\max\{\sigma(i_{\sigma}), \sigma(j_{\sigma})\})\in V(\MP_{k-1}).
 \]

 \begin{tabular}{lr}
 \small
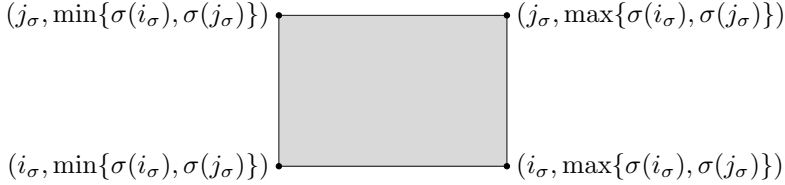
		\begin{tikzpicture}[scale=1]
		\draw[] 
        (1,1)--(4,1)--(4,3)--(1,3)--(1,1) 
        ;
        \fill[fill=gray, fill opacity=0.3] (1,1)--(4,1)--(4,3)--(1,3)--(1,1) 
        ;

        \filldraw[] (1,1) circle (1pt) node[anchor=east]  {$(i_\sigma,\min\{\sigma(i_{\sigma}), \sigma(j_{\sigma})\})$};
        \filldraw[] (4,1) circle (1pt) node[anchor=west]  {$(i_\sigma,\max\{\sigma(i_{\sigma}), \sigma(j_{\sigma})\})$};
        \filldraw[] (1,3) circle (1pt) node[anchor=east]  {$(j_\sigma,\min\{\sigma(i_{\sigma}), \sigma(j_{\sigma})\})$};
        \filldraw[] (4,3) circle (1pt) node[anchor=west]  {$(j_\sigma,\max\{\sigma(i_{\sigma}), \sigma(j_{\sigma})\})$};
		\end{tikzpicture}
 \end{tabular}  

 When $(i_{\sigma},j_{\sigma})$ satisfies~\ref{chi_a}, we get $i_{\sigma}< j_{\sigma}=l$ and $ j_{\sigma} + \max\{\sigma(i_{\sigma}), \sigma(j_{\sigma})\} \leq k + 2$. Thus, 
 \[
 (j_\sigma,\max\{\sigma(i_{\sigma}), \sigma(j_{\sigma})\})\in V(\MP_{k-1}\cup \{C_{k,l}\})
 \]
 Furthermore, since $i_{\sigma}< j_{\sigma}$ and the way $\MP_{k-1}$ is defined, we get 
  \[
 (i_\sigma,\min\{\sigma(i_{\sigma}), \sigma(j_{\sigma})\}), (i_\sigma,\max\{\sigma(i_{\sigma}), \sigma(j_{\sigma})\}), (j_\sigma,\min\{\sigma(i_{\sigma}), \sigma(j_{\sigma})\}) \in V(\MP_{k-1}\cup \{C_{k,l}\}).
 \]

 Thus, in both cases, the vertices 
 \[(i_\sigma,\sigma(i_{\sigma})), (i_\sigma,\sigma(j_{\sigma})), (j_\sigma,\sigma(i_{\sigma})), (j_\sigma,\sigma(j_{\sigma}))\in V(\MP_{k-1}\cup \{C_{k,l}\}).
 \]  
 Since $\MP_{k-1}\cup \{C_{k,l}\}$ is convex, the interval formed by these vertices is an inner interval of $\MP_{k-1}\cup \{C_{k,l}\}$.
 Therefore, $x_{{i_\sigma},\sigma(i_{\sigma})} x_{{j_\sigma},\sigma(j_{\sigma})} - x_{{i_\sigma},\sigma(j_{\sigma})} x_{{j_\sigma},\sigma(i_{\sigma})}\in I_{\MP_{k-1}\cup \{C_{k,l}\}}$.
 Hence the proof.
\end{proof}

  \begin{figure}[h]
  \small
		\centering
		\begin{tikzpicture}[scale=1]
		\draw[] 
        (1,1)--(4,1) (5,1)--(7,1) 
        (1,2)--(4,2) (5,2)--(7,2) 
        (1,3)--(4,3) (5,3)--(6,3)

        (1,4)--(4,4)
        (1,5)--(4,5)
        (1,6)--(3,6)
        (1,7)--(2,7)
       
        (1,1)--(1,3) (1,4)--(1,7)
        (2,1)--(2,3) (2,4)--(2,7)
        (3,1)--(3,3) (3,4)--(3,6)
        (4,1)--(4,3) (4,4)--(4,5)
        
        (5,1)--(5,3)
        (6,1)--(6,3)
        (7,1)--(7,2)

        ;
        \draw[dotted, thick]
        (4,1)--(5,1)
        (4,2)--(5,2)
        (4,3)--(5,3)
        (4,4)--(5,4)

        (1,4)--(1,3)
        (2,4)--(2,3)
        (3,4)--(3,3)
        (4,4)--(4,3)
        (5,3)--(5,4)

        ;

        \filldraw[] (1.5,6.5) circle (0pt) node[anchor=center]  {$C_{k,1}$};
        \filldraw[] (2.5,5.5) circle (0pt) node[anchor=center]  {$C_{k,2}$};
        \filldraw[] (3.5,4.5) circle (0pt) node[anchor=center]  {$C_{k,3}$};
        \filldraw[] (6.5,1.5) circle (0pt) node[anchor=center]  {$C_{k,n_k-1}$};

        \filldraw[] (1,1) circle (1pt) node[anchor=north]  {$(i_k+1,j_k+1)$};
        \filldraw[] (7,1) circle (1pt) node[anchor=north]  {$v_{k,n_k} = (i_k+n_k, j_k+1)$};
        \filldraw[] (1,7) circle (1pt) node[anchor=east]  {$v_{k,1}=(i_k+1, j_k + n_k)$};
        \filldraw[] (1,6) circle (1pt) node[anchor=east]  {$(i_k+1, j_k + n_k-1)$};

		\end{tikzpicture}
		\caption{}%
	\label{fig:C_knbd}
\end{figure}

\begin{lemma}\label{lem:subpolyomino}
    Let $\MP$ be a ladder polyomino. Under the notations of Discussion~\ref{setup:ladder}, for each $1\leq k\leq r$, 
    \begin{asparaenum}
    \item 
        if $\MC_k\neq \varnothing$, then $\MP_k$ contains the polyomino shown in Figure~\ref{fig:C_knbd} as a sub-polyomino.
    \item \label{lem:Pkprime}
        if $\MP_k\neq \varnothing$, $\MP_k$ is a simple sub-polyomino of $\MP.$
    \end{asparaenum}
\end{lemma}

\begin{proof}
    $(1)$ Assuming the notations from Discussion~\ref{setup:ladder},
    we note that if $\MC_k\neq \varnothing$, then for all $1\leq l\leq n_k-1$, the cell $C_{k,l}$ has upper-left vertex $(i_k+l, j_k + n_k+1-l)$.
    Given the definition of $\MP_k$, it suffice to show that the polyomino shown in Figure~\ref{fig:C_knbd} is a sub-polyomino of $\MP$. To prove that, we need to show that all cells with upper-left vertex $(i_k+u, j_k+v)$ are in $\MP$, where $1\leq u\leq n_k-1, 2\leq v\leq n_k,$ and $u+v\leq n_k+1.$ 
    We proceed by downward induction on $v$.
    When $v=n_k,$ we have $u=1$. In this case, the cell with upper-left vertex $(i_k+1, j_k+n_k)$ is $C_{k,1}$ which is in $\MP.$
    Let $v<n_k$ and by induction hypothesis, assume that all cells with upper-left vertex $(i_k+u,j_k+v+1)$ are in $\MP$, where $1\leq u\leq n_k-1$ and $u+v+1\leq n_k+1.$
    In particular, the cell with upper-left vertex $(i_k+1,j_k+v+1)$ is in $\MP$, implying that the vertex $(i_k+1,j_k+v)\in V(\MP).$
    
    Let $a_l$ be the left-most vertex of $\MP$ such that $a_{l,2}=j_k+v$. Given that $(i_k+1,j_k+v)\in V(\MP),$ it follows that that $a_{l,1}\leq i_k+1.$
    Since $\MP$ is a ladder polyomino, the cell with upper-left vertex $a_l=(a_{l,1},j_k+v)$ is in $\MP$.
    Moreover, since $\MP$ is horizontally convex and the cell $C_{k,n_k-v+1}$ with upper-left vertex $(i_k+n_k-v+1,j_k+v)$ is also in $\MP$, we get that for all $m$, $a_{l,1}\leq m\leq i_k+n_k-v+1,$ the cell with upper-left vertex $(m,j_k+v)$ is in $\MP$. Hence the proof as $a_{l,1}\leq i_k+1.$

    $(2)$ Note that $\MP_1=\varnothing$, and $\MP_2$ is a cell.
    Let $k\geq 3$. We use induction, assuming that $\MP_{k-1}$ is a simple sub-polyomino of $\MP.$
    By $(1)$, the cells of $\MC_k$ are attached to $\MP_{k-1}$ in such a way that $\MP_k$ is a polyomino.
    Thus, $\MP_k$ is a sub-polyomino of $\MP.$ Also, given our definition of the set $V_i$'s, $\MP_k$ is simple. 
\end{proof}
 
We are now ready to prove our main theorem of this section.

\begin{proof}[Proof of Theorem~\ref{thm:ladderpolyomino}]
We show that $I_\MP\in C_f$, where $f$ is as defined in Discussion~\ref{setup:ladder}.
By Lemma \ref{lem:one_sided_initial}, $\ini_{<}(f) = \prod_{v \in V(\MP)} x_v.$
Assume the notations of Discussion~\ref{setup:ladder}.
Note that $V_1=\{(i_1+1, j_1 + 1)\}$ is singleton, and $\MP_1 =\varnothing$.
We proceed by induction on $k$ for $k\geq2$ to show that $I_{\MP_k}\in C_f$.
Since $\MP$ is a ladder polyomino, $V_2=\{(i_1+1, j_1 + 2),(i_1+2, j_1 + 1)\},$ and the polyomino ideal $I_{\MP_2}$ is generated by $x_{i_1+1, j_1 + 1}x_{i_1+2, j_1 + 2}-x_{i_1+1, j_1 + 2}x_{i_1+2, j_1 + 1}$ which is $-f_2$. 
Thus, $I_{\MP_2}\in C_f.$

Now, let $k>2$ and by induction hypothesis, assume that $I_{\MP_{k-1}}\in C_f$. Recall that $\MP_k = \MP_{k-1}\cup \MC_k$, where $\MC_k= \{C_{k,1}, C_{k,2}, \ldots, C_{k,{n_k - 1}}\}$. If $\MC_k=\varnothing,$ then $\MP_k = \MP_{k-1}$ and we are done. So we may assume that $\MC_k\neq \varnothing$.
We first show that 
\[I_{\MP_k} = \sum_{C\in \MC_k} I_{\MP_{k-1} \cup \{C\}}.\]
Clearly, $\sum_{C\in \MC_k} I_{\MP_{k-1} \cup \{C\}} \subseteq I_{\MP_k}.$
Given our definition of the set $V_i$'s, we observe that $\MP_k$ does not contain a cell with lower-left vertex $v_{k,i}$ for all $2\leq i\leq n_k.$
Consequently, there are no inner intervals in $\MP_k$ that contain at least two cells from the set $\MC_k$. 
Thus, for each inner interval $I$ in $\MP_k$, there exists an $C\in \MC_k$ such that $I$ is an inner interval of $\MP_{k-1} \cup \{C\}$.
Thus, $I_{\MP_k} \subseteq \sum_{C\in \MC_k} I_{\MP_{k-1} \cup \{C\}}$.
Now, it suffices to show that $I_{\MP_{k-1} \cup \{C\}} \in C_f$ for all $C\in \MC$. 

We claim that $f_{k} \in I_{\MP_{k-1} \cup \{C\}}$ for $C\in \MC_k$.
Recall  from the proof of Lemma~\ref{lem:one_sided_initial},
\[
	f_{k} = \sum_{\sigma \in S_{n_k}} (-1)^{\sgn(\sigma)} x_{i_k+1, j_k + \sigma(1)} x_{i_k+2, j_k + \sigma(2)} \cdots x_{i_k+n_k, j_k + \sigma(n_k)}.
 \]
 where $S_{n_k}$ denotes the symmetric group on $n_k$ symbols, and $\sgn(\sigma)$ denotes the sign of the permutation $\sigma$, which is $1$ if $\sigma$ is even and $-1$ if $\sigma$ is odd. 
 Some of the monomials in the summation might not be present.
By Lemma~\ref{lem:subpolyomino}, $P_k$ contains the polyomino shown in Figure~\ref{fig:C_knbd} as a sub-polyomino, call it $\MQ$.
Moreover, by the proof of Lemma~\ref{lem:subpolyomino}, the polyomino $\MP_{k-1}$ contains the polyomino $\MQ\setminus \MC_k$ as a sub-polyomino.
So by Lemma~\ref{lem:detf_k}, $f_{k} \in I_{\MP_{k-1} \cup \{C\}}$ for all $C\in \MC_k$.

Observe that $I_{\MP_{k-1}}$ and $I_{\MP_{k-1} \cup \{C\}}$ are prime ideals for every $C\in \MC_k$ by Lemma~\ref{lem:subpolyomino} and~\cite[Corollary\ 2.2]{HM14}.
Also, note that $(f_{k}) \in C_f$. 
Therefore, $I_{\MP_{k-1}} + (f_{k}) \in C_f$. 
Since $f_{k} \notin I_{\MP_{k-1}}$ by Lemma~\ref{lem:detf_k}, we have $\height(I_{\MP_{k-1}} + (f_{k})) = \height(I_{\MP_{k-1}}) + 1$.
Thus, $I_{\MP_{k-1} \cup \{C\}}$ is a minimal prime over $I_{\MP_{k-1}} + (f_{k})$, 
which implies $I_{\MP_{k-1} \cup \{C\}} \in C_f$.
Hence the proof.
\end{proof}

\begin{Remark}\label{rem:differentpolynomial}
    Under the notations of Discussion~\ref{setup:ladder}, define 
    \[
    g := \prod_{\substack{k \in \{1,\ldots,s\}\\
                 \text{with} \  |V_k|\geq 2}} f_{k}.
    \]
    Then, one can use similar argument of proof of Theorem~\ref{thm:ladderpolyomino} to get that $I_\MP\in C_g$ if $\MP$ is a ladder polyomino. It holds because if $|V_k|=1$ for some $k,$ then $\MP_k = \MP_{k-1}.$
 \end{Remark}

Given a cell $C = [a,a+(1,1)]$, we define $\sp(C) = \{x\in \reals^2: a\leq x\leq a+(1,1)\}$.  Let $\MP$ be a collection of cells, define $\sp(\MP) = \cup_{C\in \MP} \sp(C)\subset\reals^2$.
Let $\MP$ be a ladder polyomino with the set of left-most vertices $\{a_1=(a_{1,1},a_{1,2}),\ldots,a_r=(a_{r,1},a_{r,2})\}$. 
Recall that $a_{i+1,2} = a_{i,2}+1$ for all $1\leq i \leq r-1$. Let $b = \max\{v\in \naturals\mid (v,a_{1,2})\in V(\MP)\}.$ Since $\MP$ is horizontally convex, the horizontal edge interval $[a_1,(b,a_{1,2})]$ is in $\MP.$ We say that $\MP$ is \define{based on} a polyomino $\MQ,$ if $\sp(\MP)\cap \sp(\MQ) = [a_1,(b,a_{1,2})]$, here $[a_1,(b,a_{1,2})]$ is the closed interval in $\reals^2$.

In the proof of Theorem~\ref{thm:ladderpolyomino}, we utilized the hypothesis that $\MP$ is a ladder polyomino to apply Lemma~\ref{lem:subpolyomino}. The proof goes through for any polyomino for which Lemma~\ref{lem:subpolyomino} holds true. Notably, instead of~\eqref{lem:Pkprime} of Lemma~\ref{lem:subpolyomino}, it suffices to assume that the ideals $I_{\MP_k}$ and $I_{\MP_k\cup \{C\}}$ are prime for all $C\in \MC_{k+1}.$ This leads us to the following result:

  \begin{figure}[h]%
		\centering
		\begin{tikzpicture}[scale=1]
		\draw[] 
        (0,1.5)--(11.75,1.5)--(11.75,2)--(12.5,2)--(12.5,2.5)--(13,2.5)--(13,3.2)

        (1,3) -- (4.5,3)--(4.5,3.1)--(9,3.1)--(9,3.2)--(13,3.2)
        (1,3)--(1,2.5)--(.6,2.5)--(.6,2)--(0,2)--(0,1.5)

        (1.5,3)--(1.5,4)--(2,4)--(2,5)--(2.5,5)--(2.5,6)--(3.5,6)--(3.5,5.5)--(4,5.5)--(4,4.5)--(3.5,4.5)--(3.5,3.5)--(4,3.5)--(4,3)
        
        (5,3.1)--(5,3.5)--(5.5,3.5)--(5.5,5)--(6,5)--(6,6)--(7,6)--(7,5.5)--(7,5)--(7,4.5)--(7.5,4.5)--(7.5,4)--(7,4)--(7,3.5)--(7.5,3.5)--(7.5,3.1)

        (10,3.2)--(10,4)--(10.5,4)--(10.5,5)--(11,5)--(11,5.5)--(11.5,5.5)--(11.5,6)--(12,6)--(12,5.5)--(12.5,5.5)--(12.5,5)--(12,5)--(12,4.5)--(12.5,4.5)--(12.5,4)--(12,4)--(12,3.5)--(12.5,3.5)--(12.5,3.2)

   ;

       \filldraw[] (7,2.25) circle (0pt) node[anchor=center]  {$\MI$};
       \filldraw[] (2.75,4) circle (0pt) node[anchor=center]  {$\MR_1$};
       \filldraw[] (6.25,4) circle (0pt) node[anchor=center]  {$\MR_2$};
       \filldraw[] (9,4) circle (0pt) node[anchor=center]  {$\cdots$};
       \filldraw[] (11.25,4) circle (0pt) node[anchor=center]  {$\MR_r$};

		\end{tikzpicture}
		\caption{}%
	\label{fig:ladders}
\end{figure}

\begin{proposition}\label{prop:ladderson interval}
Let $\MP$ be the polyomino as shown in Figure~\ref{fig:ladders}, where $\MI$ is a parallelogram polyomino and $\MR_1,\ldots,\MR_r$ are ladder sub-polyominoes based on $\MI$ such that $V(\MR_i)\cap V(\MR_j)=\varnothing$ for $i\neq j$.
Then, $I_{\MP} \in C_f$, where $f$ is as defined in Discussion~\ref{setup:ladder}. In fact $I_{\MP} \in C_g$, where $g$ is as defined in Remark~\ref{rem:differentpolynomial}.
\end{proposition}

\begin{proof}
Follow Discussion~\ref{setup:ladder} to define $V_k$'s and $\MP_k$'s for $\MP$.  For $1\leq k\leq s,$ note that the vertex set $V_k$ is contained in $\MI\cup \MR_i$ for some $1\leq i\leq r$. 
If $V_k$ is contained in $V(\MI)$ (resp. $V(\MR_i)$) then by Lemma~\ref{lem:subpolyomino}, $\MP_k$ contains the polyomino shown in Figure~\ref{fig:C_knbd} as a sub-polyomino because $\MI$ (resp. $\MR_i$) is a ladder polyomino.
If $V_k \cap V(\MR_i)\neq \varnothing$ and $V_k \cap V(\MI)\neq \varnothing$, then define a sub-polyomino $\MI_k$ of $\MP$ as 
the cells of $\MP_k$ contained in $\MI$ or $\MR_i.$
Note that $\MI_k$ is a ladder polyomino because $\MI$ and $\MR_i$ are ladder polyominoes, and $\MR_i$ is based on $\MI$. So by Lemma~\ref{lem:subpolyomino}, $\MI_k$ contains the polyomino shown in Figure~\ref{fig:C_knbd} as a sub-polyomino. 
By following the argument of Lemma~\ref{lem:subpolyomino}, and the assumption $V(\MR_i)\cap V(\MR_j)=\varnothing$ for $i\neq j$, we get that $\MP_k$ is a simple sub-polyomino of $\MP.$
We can now use the argument of proof of Theorem~\ref{thm:ladderpolyomino} to show that $I_{\MP} \in C_f$.
Moreover, $I_{\MP} \in C_g$ by Remark~\ref{rem:differentpolynomial}.
\end{proof}

\begin{proposition}\label{prop:groebner}
    Let $\MP$ be as defined in Proposition~\ref{prop:ladderson interval} such that $\MR_1,\ldots,\MR_r$ are parallelogram polyominoes. Then, the set of all quadratic binomials corresponding to inner intervals of $\MP$ form a Gr\"{o}bner basis of $I_{\MP}$ under the order $<$ as defined in Discussion~\ref{setup:ladder}. 
\end{proposition}
\begin{proof}
    We apply Buchberger’s criterion. Let $f = x_ax_d-x_bx_c$ and $g = x_{a'}x_{d'}-x_{b'}x_{c'}$ be two quadratic binomials corresponding to inner intervals $I$ and $J$ of $\MP$ respectively. 
    Let $a$ and $d$ (resp. $a'$ and $d'$) be the anti-diagonal vertices of $I$ (resp. $J$).
    Then, $\ini_<(f) = x_ax_d$ and $\ini_<(g) = x_{a'}x_{d'}$. If $\{a,d\}\cap \{a',d'\} =\varnothing,$ then by~\cite[Lemma\ 1.27]{HHO18}, the $S$-pair $S(f,g)$ reduces to zero. On the other hand, if $\{a,d\} =  \{a',d'\}$, then the inner interval $I=J$; thus $f=g$. So, we may assume that $|\{a,d\}\cap \{a',d'\}| =1,$ say $a=a'.$ 

    Since $\MI, \MR_1,\ldots,\MR_r$ are parallelogram polyominoes, the $K$-algebras $K[\MI], K[\MR_1],\ldots,K[\MR_r]$ are the Hibi rings associated to the vertex sets $V(\MI), V(\MR_1),\ldots,V(\MR_r)$ respectively, which are distributive lattices (see~\cite{KV21}, \cite[Remark\ 2.4]{QRR22}).
    When both inner intervals $I$ and $J$ are contained in one of the sub-polyominoes $\MI, \MR_1,\ldots,\MR_r.$ Then, $S(f,g)$ reduces to zero by~\cite[Theorem\ 6.17]{HHO18} (after noting that the monomial order $<$ is a ``compatible order" as defined in~\cite[Section\ 6.2]{HHO18}).

    Consider the case when these inner intervals $I$ and $J$ are not contained in one of these sub-polyominoes. 
    If $\sp(I) \cap \sp(J)\neq \{a\}$, then $I\cup J$ is a parallelogram polyomino. Again, by~\cite[Theorem\ 6.17]{HHO18},  $S(f,g)$ reduces to zero. So, we may assume that $\sp(I) \cap \sp(J) = \{a\}$.
    Since $I$ and $J$ intersects, these intervals must be contained in the sub-polyomino consists of cells of $\MI$ and $\MR_i$ for some $1\leq i\leq r.$ 
    We may assume that the inner interval $I$ intersects with $\MI$ and the inner interval $J$ intersects with $\MR_i$.
    Note that $d\notin V(\MR_i)$; otherwise, both $I$ and $J$ must be contained in $\MR_i$ because $\MR_i$ is based on $\MI$. Similarly, $d'\notin V(\MI)$. 
    Consequently, $b\notin V(\MR_i)$ and $b'\notin V(\MI)$ because $\MR_i$ is based on $\MI$.
    Without loss of generality, assume that $b$ (resp. $b'$) is the lower-left vertex of $I$ (resp. $J$). 

    If $b'\in V(\MI)$, then $b'\wedge b \in V(\MI)$ because $V(\MI)$ is a sub-lattice of $\naturals^2,$ where $\wedge$ denotes the meet in the distributive lattice $\naturals^2$. Since $\MI$ is convex, the interval $[b'\wedge b, a]$ is an inner interval of $\MI$. Thus, $[b'\wedge b, c']$ and $[b'\wedge b, c]$ are inner intervals of $\MI\cup \MR_i$. So we get that,
    \[
    S(f,g) = x_{c'}(x_{d}x_{b'}-x_{c} x_{b'\wedge b}) - x_c(x_{b} x_{d'}-x_{c'} x_{b'\wedge b})
    \]
    which reduced to zero.
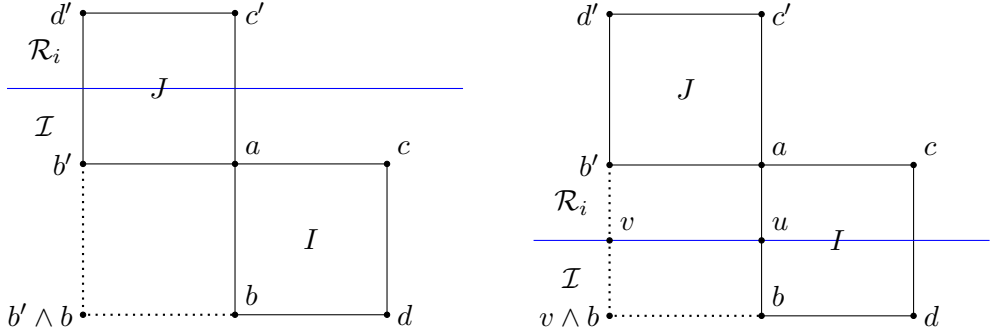
\[
\begin{tabular}{lr}
      \begin{tikzpicture}[scale = 1]
      \draw[]
            (3,3)--(5,3)--(5,5)--(3,5)--(3,3)
            (5,3)--(5,1)--(7,1)--(7,3)--(5,3)
            ;

       \draw[thick, dotted] (3,3)--(3,1)--(5,1);

       \draw[blue] (2,4)--(8,4);

       \filldraw[] (2.5,4.5) circle (0pt) node[anchor=center]  {$\MR_i$};
       \filldraw[] (2.5,3.5) circle (0pt) node[anchor=center]  {$\MI$};
       
        \filldraw[] (6,2) circle (0pt) node[anchor=center]  {$I$};
       \filldraw[] (4,4) circle (0pt) node[anchor=center]  {$J$};       

       \filldraw[] (5,3) circle (1pt) node[anchor=south west]  {$a$};
       \filldraw[] (7,1) circle (1pt) node[anchor=west]  {$d$};
       \filldraw[] (7,3) circle (1pt) node[anchor=south west]  {$c$};
       \filldraw[] (5,1) circle (1pt) node[anchor=south west]  {$b$};
       
       \filldraw[] (3,5) circle (1pt) node[anchor=east]  {$d'$};
       \filldraw[] (5,5) circle (1pt) node[anchor=west]  {$c'$};
       \filldraw[] (3,3) circle (1pt) node[anchor=east]  {$b'$};
       
       \filldraw[] (3,1) circle (1pt) node[anchor=east]  {$b'\wedge b$};
       \end{tikzpicture}
       \quad & \quad
    \begin{tikzpicture}[scale = 1]
      \draw[]
            (3,3)--(5,3)--(5,5)--(3,5)--(3,3)
            (5,3)--(5,1)--(7,1)--(7,3)--(5,3)
            ;

       \draw[thick, dotted] (3,3)--(3,1)--(5,1);

       \draw[blue] (2,2)--(8,2);

       \filldraw[] (2.5,2.5) circle (0pt) node[anchor=center]  {$\MR_i$};
       \filldraw[] (2.5,1.5) circle (0pt) node[anchor=center]  {$\MI$};

        \filldraw[] (6,2) circle (0pt) node[anchor=center]  {$I$};
       \filldraw[] (4,4) circle (0pt) node[anchor=center]  {$J$};

       \filldraw[] (5,2) circle (1pt) node[anchor=south west]  {$u$};
       \filldraw[] (3,2) circle (1pt) node[anchor=south west]  {$v$};
       \filldraw[] (3,1) circle (1pt) node[anchor=east]  {$v\wedge b$};

       \filldraw[] (5,3) circle (1pt) node[anchor=south west]  {$a$};
       \filldraw[] (7,1) circle (1pt) node[anchor=west]  {$d$};
       \filldraw[] (7,3) circle (1pt) node[anchor=south west]  {$c$};
       \filldraw[] (5,1) circle (1pt) node[anchor=south west]  {$b$};
       
       \filldraw[] (3,5) circle (1pt) node[anchor=east]  {$d'$};
       \filldraw[] (5,5) circle (1pt) node[anchor=west]  {$c'$};
       \filldraw[] (3,3) circle (1pt) node[anchor=east]  {$b'$};
       \end{tikzpicture}      
\end{tabular}
\]
Consider the case when $b'\notin V(\MI)$. So, $b'\in V(\MR_i)$. Since $\MR_i$ is based on $\MI$, we get that $a\in V(\MR_i)$.
Since $a\in V(\MR_i),$ $b\in V(\MI)$ and the inner interval $I$ is in $\MI\cup\MR_i,$ we get that there is a $u=(u_1,u_2)\in V(\MI)\cap V(\MR_i)$ such that $b\leq u\leq a.$ Since $V(\MR_i)$ is a sub-lattice of $\naturals^2$, we get that $u\wedge b'\in V(\MR_i)$, call it $v=(v_1,v_2)$. Since $\MR_i$ is convex, the interval $[v,a]$ is an inner interval of $\MR_i.$ Moreover, $v\in V(\MI)$ because $u_2 = v_2$ and $\MR_i$ is based on  $\MI$. Thus, $v\wedge b \in V(\MI)$. Furthermore, the interval $[v\wedge b, u]$ is an inner interval of $\MI$ because $\MI$ is convex. Thus, $[v\wedge b, c']$ and $[v\wedge b, c]$ are inner intervals of $\MI\cup \MR_i.$ 
So we get that,
    \[
    S(f,g) = x_{c'}(x_{d}x_{b'}-x_{c} x_{v\wedge b}) - x_c(x_{b} x_{d'}-x_{c'} x_{v\wedge b})
    \]
    which reduced to zero.
\end{proof}

\section{Knutson thin polyominoes II}\label{sec:thin2}

In this section, we apply the tools developed in the previous section to prove that binomials ideals of a class of thin collection to cells are Knutson. We also compute its Gr\"{o}bner basis. As a consequence we prove that when this collection of cells is a polyomino, then the corresponding polyomino ideal is prime. 

\begin{theorem}\label{thm:thingap2}
  Let $\MP$ be a thin collection of cells such that 
   \begin{asparaenum}\label{thm:thingap2:3cells}
       \item
           it does not contain cells as in Figure~\ref{fig:P_k4elements};

       \item \label{thm:thingap2:vertical}
            if $I$ is a maximal inner interval of $\MP$ with vertices $(i_1,j_1)$, $(i_1+1,j_1 )$, $(i_1,j_2 )$ and $(i_1+1, j_2 )$ (see Figure~\ref{fig:verticalinterval}), then the cell with lower-right vertex $(i_1,j_2)$ and the cell with upper-left vertex $(i_1+1,j_1)$ are not contained in $\MP$;

        \item \label{thm:thingap2:horizontal}
            if $I$ is a maximal inner interval of $\MP$ with vertices $(i_1,j_1)$, $(i_1,j_1+1)$, $(i_2,j_1 )$ and $(i_2, j_1+1)$ (see Figure~\ref{fig:horizontalinterval}), then the cell with lower-right vertex $(i_1,j_1+1)$ and the cell with upper-left vertex $(i_2,j_1)$ are not contained in $\MP$. 
       \end{asparaenum}
  Then, there exists an $f\in S$ such that \(I_{\MP} \in C_f\), that is, \(I_{\MP}\) is a Knutson ideal.
\end{theorem}

  \begin{figure}[h]
 \begin{subfigure}[t]{7cm}
    \centering
		\begin{tikzpicture}[scale=.8]
  \draw[]
      (1,4)--(2,4)--(2,3)--(1,3)--(1,4)
        (2,3)--(3,3)--(3,2)--(2,2)--(2,3)
        (3,2)--(4,2)--(4,1)--(3,1)--(3,2);

        \filldraw[] (1,4) circle (1pt) node[anchor=east]  {$(i,j)$};
        \filldraw[] (2,3) circle (1pt) node[anchor=south west]  {$(i+1,j-1)$};
        \filldraw[] (3,2) circle (1pt) node[anchor=south west]  {$(i+2,j-2)$};
        \filldraw[] (4,1) circle (1pt) node[anchor=west]  {$(i+3,j-3)$};

	\end{tikzpicture}
    \caption{}\label{fig:P_k4elements}
    \end{subfigure}
    \quad
  \begin{subfigure}[t]{5cm}
    \centering
		\begin{tikzpicture}[scale=.8]
  \draw[]
      (1,1)--(1,4)--(2,4)--(2,1)--(1,1)
        ;

        \filldraw[] (1,1) circle (1pt) node[anchor=east]  {$(i_1,j_1)$};
        \filldraw[] (1,4) circle (1pt) node[anchor=east]  {$(i_1,j_2)$};
        \filldraw[] (2,1) circle (1pt) node[anchor=west]  {$(i_1+1,j_1)$};
        \filldraw[] (2,4) circle (1pt) node[anchor=west]  {$(i_1+1,j_2)$};
	\end{tikzpicture}
    \caption{} \label{fig:verticalinterval}
    \end{subfigure}
     \quad
  \begin{subfigure}[t]{5cm}
    \centering
		\begin{tikzpicture}[scale=.8]
  \draw[]
      (1,1)--(1,2)--(4,2)--(4,1)--(1,1)
        ;

        \filldraw[] (1,1) circle (1pt) node[anchor=east]  {$(i_1,j_1)$};
        \filldraw[] (1,2) circle (1pt) node[anchor=east]  {$(i_1,j_1+1)$};
        \filldraw[] (4,1) circle (1pt) node[anchor=west]  {$(i_2,j_1)$};
        \filldraw[] (4,2) circle (1pt) node[anchor=west]  {$(i_2,j_1+1)$};
	\end{tikzpicture}
    \caption{} \label{fig:horizontalinterval}
    \end{subfigure}
        \caption{}
 \end{figure}     

Consider the polyomino $\MP$ shown in Figure~\ref{fig:examplethin}.
$\MP$ satisfies the hypothesis of Theorem~\ref{thm:thingap2} even though 
it contain a collection of cells isomorphic to as cells shown in Figure~\ref{fig:P_k4elements}.
That is why the labeling in Figure~\ref{fig:P_k4elements} is important.
 \begin{figure}[h]
    \centering
		\begin{tikzpicture}[scale=.8]
  \draw[]
      (1,2)--(2,2)--(2,1)--(1,1)--(1,2)
        (2,3)--(3,3)--(3,2)--(2,2)--(2,3)
        (3,4)--(4,4)--(4,3)--(3,3)--(3,4)
        (2,3)--(2,4)--(3,4)
        (1,2)--(1,3)--(2,3)
        ;

	\end{tikzpicture}
        \caption{}\label{fig:examplethin}
 \end{figure}  

\begin{proof}
For $\MP$, let $f$ be as defined in Discussion~\ref{setup:ladder}. 
We show that $I_\MP \in C_f$.
By Lemma \ref{lem:one_sided_initial}, $\ini_{<}(f) = \prod_{v \in V(\MP)} x_v.$
Let  $\{J_1, \ldots, J_r\}$ be the set of all maximal inner intervals of $\MP$.
Then $I_{\MP} = \sum_{i=1}^{r} I_{J_i}$, where $I_{J_i}$ denotes the polyomino ideal associated to the maximal inner interval $J_i$. 
So, it suffice to show that each $I_{J_i} \in C_f$ for all $1\leq i\leq r$. 
Let $J$ be a maximal inner interval of $\MP$. 
Let $\MC_{i_1},\ldots, \MC_{i_m}$ be the set of cells (as defined in Discussion~\ref{setup:ladder}) such that $\MC_{i_k}\cap J \neq \varnothing$ for all $k$. 
Since $\MP$ is thin, $\MC_{i_k}\cap J$ is a cell, say $C_k$ for all $1\leq k\leq m.$ 
Thus, $J$ is a maximal inner interval consisting of cells $C_1,\ldots,C_m$.
Without loss of generality, assume that $i_1<i_2<\cdots<i_m.$ Thus, $C_1$ is an end-cell of $J$, and $C_k$ and $C_{k+1}$ are neighbour cells for all $1\leq k \leq m-1$.

Fix a $k$ with $1\leq k\leq m$. 
Since $\MC_{i_k}\neq \varnothing,$ we have $|V_{i_k}|\geq 2.$ 
Also, since $\MP$ does not contain cells as in Figure~\ref{fig:P_k4elements}, $|V_{i_k}|\leq  3.$
If $|V_{i_k}| = 2$, then it follows that $f_{i_k}= I_{C_k}\subseteq I_J$, where $I_{C_k}$ denotes the polyomino ideal associated to the cell $C_k$, and $f_{i_k}$ is as defined in Discussion~\ref{setup:ladder}. 
Now, consider the case when $|V_{i_k}| = 3$ and $1<k<m$. Since $|V_{i_k}| = 3$, it follows that $|\MC_{i_k}| =2$; we know $C_k\in \MC_{i_k}$, let the other cell in $\MC_{i_k}$ is denoted by $D$. The cell $D$ is either a neighbour cell of $C_{k-1}$ or a neighbour cell of $C_{k+1}.$  
If $D$ is a neighbour cell of $C_{k-1}$ (see Figure~\ref{fig:proofthinfirst}), we have $f_{i_k} \in I_{C_{k-1} \cup C_k}$ by Lemma~\ref{lem:detf_k}. If $D$ is a neighbour cell of $C_{k+1}$ (see Figure~\ref{fig:proofthinsecond}), then again by Lemma~\ref{lem:detf_k}, we get that $f_{i_k} \in I_{C_k \cup C_{k+1}}$. Thus, $f_{i_k} \in I_J$ in both cases.
 
		\begin{figure}[h]
			\centering
			\begin{subfigure}[t]{0.45\textwidth}
				\centering
				\begin{tikzpicture}[scale=1]
					\draw[] (0.0,0.0) -- (0,2);
					\draw[] (0,0) -- (2,0);
					\draw[] (1,0) -- (1,2);
					\draw[] (0,1) -- (2,1);
					\draw[] (0,2) -- (1,2);
					\draw[] (2,0) -- (2,1);

                    \draw[thick, dotted]
                    (2,0)--(3,0)--(3,1)--(2,1)
                    (0,0)--(-1,0)--(-1,1)--(0,1)

                    ;
                    \draw[color=blue] (0,2) -- (2,0);
                    
                    \filldraw[] (3.5,.5) circle (0pt) node[anchor=center]  {$J$};
                    \filldraw[] (0,2) circle (0pt) node[anchor=south east]  {$V_{i_k}$};
					\node at (0.5,1.5) {$D$};
                    \node at (0.5,0.5) {$C_{k-1}$};
					\node at (1.5,0.5) {$C_{k}$};
					
				\end{tikzpicture}
				\caption{}\label{fig:proofthinfirst}
			\end{subfigure}
			\quad
			\begin{subfigure}[t]{0.45\textwidth}
				\centering
				\begin{tikzpicture}[scale=1]
					\draw[] (1,0) -- (2,0) -- (2,2) -- (0,2) -- (0,1) -- (1,1) -- (1,0);
					\draw[] (1,2) -- (1,1) -- (2,1);

                    \draw[thick, dotted]
                    (2,1)--(3,1)--(3,2)--(2,2)
                    (0,1)--(-1,1)--(-1,2)--(0,2)
                    ;
                    \draw[color=blue] (0,2) -- (2,0);

                    \filldraw[] (3.5,1.5) circle (0pt) node[anchor=center]  {$J$};
                    \filldraw[] (0,2) circle (0pt) node[anchor=south east]  {$V_{i_k}$};
					\node at (1.5,1.5) {$C_{k+1}$};
					\node at (0.5,1.5) {$C_{k}$};
                    \node at (1.5,0.5) {$D$};
					
				\end{tikzpicture}
				\caption{}\label{fig:proofthinsecond}
			\end{subfigure}
			\caption{}\label{fig:proofthin}
		\end{figure}

Now consider the case when $|V_{i_k}| = 3$ and $k\in \{1,m\}$. we know $C_1\in \MC_{i_1}$ (resp. $C_m\in \MC_{i_m}$), let the other cell in $\MC_{i_1}$ (resp. $\MC_{i_m}$) is denoted by $D$.  
We know $C_1\in \MC_{i_1}$, let the other cell in $\MC_{i_1}$ is denoted by $D$. 
When $J$ is as defined in \eqref{thm:thingap2:vertical} of theorem. 
Since $\MP$ does not contain the cell with lower-right vertex $(i_1,j_2)$, it follows that $V_{i_1} = \{(i_1,j_2),(i_1+1,j_2-1),(i_1+2,j_2-1)\}$. Thus, $D$ is a neighbour cell of $C_2$. 
When $J$ is as defined in \eqref{thm:thingap2:horizontal} of theorem. 
Since $\MP$ does not contain the cell with lower-right vertex $(i_1,j_1+1)$, we get that $V_{i_1} = \{(i_1,j_1+1),(i_1+1,j_1),(i_1+2,j_1-1)\}$. Thus, $D$ is a neighbour cell of $C_2$.
Similarly, the other cell of $\MC_{i_k}$ is a neighbour cell of $C_{k-1}$ in both subcases.
Using the earlier argument, we get that $f_{i_1}, f_{i_m} \in I_J$. Thus, $(f_{i_1}, f_{i_2}, \ldots, f_{i_m})\subseteq I_{J}$. 

Note that $(f_{i_k}) \in C_f$ for all $1 \leq k \leq m$, and consequently the ideal $(f_{i_1}, f_{i_2}, \ldots, f_{i_m}) \in C_f$. 
Since $\ini_{<}(f_{i_k})=\prod_{v\in V_{i_k}}x_v$ for all $k$ and the set $V_{i_k}$'s are disjoint, we get that $\{f_{i_1}, f_{i_2}, \ldots, f_{i_m}\}$ forms a Gr\"{o}bner basis of $(f_{i_1}, f_{i_2}, \ldots, f_{i_m})$ by ~\cite[Lemma\ 1.27]{HHO18}, and also we get the that height of the ideal $(\ini_{<}(f_{i_1}), \ldots, \ini_{<}(f_{i_m})$ is $m$;
Thus $\height((f_{i_1}, f_{i_2}, \ldots, f_{i_m})) = m$.
Since $I_J$ is a prime polyomino ideal of height $m$, it is a minimal prime ideal of $(f_{i_1}, f_{i_2}, \ldots, f_{i_m})$;
therefore, $I_J\in C_f$.
Hence the proof.
\end{proof}

 \begin{corollary}\label{coro:thinprime}
	Let $\MP$ be a thin collection of cells as defined in Theorem~\ref{thm:thingap2}. Then, the set of all binomials corresponding to inner intervals of $\MP$ forms a Gr\"{o}bner basis of $I_{\MP}$ under the order $<$ as defined in Discussion~\ref{setup:ladder}. 
    Moreover, when $\MP$ is a polyomino, then $I_{\MP}$ is in fact a prime ideal.
\end{corollary}

\begin{proof}
Recall from the proof of Theorem \ref{thm:thingap2}, that $I_{\MP}=\sum I_J$ where the summation runs over all the maximal inner intervals $J$ of $\MP$, and each $I_J\in C_f$ where $f$ is as defined in Discussion~\ref{setup:ladder}. Therefore, the union of Gr\"obner bases of Knutson ideals $I_J$ under the monomial order $<$, is a Gr\"obner basis of $I_\MP$ under the same monomial order. So it suffice to show that the set of all binomials corresponding to inner intervals of a maximal inner interval $J$ of $\MP$ form a Gr\"{o}bner basis of $I_J$ under the order $<$.

When $J$ is an inner interval, the $K$-algebra $K[J]$ is the Hibi ring associated to the vertex set $V(J)$ which is a distributive lattice (see~\cite{KV21}, \cite[Remark\ 2.4]{QRR22}). By~\cite[Theorem\ 6.17]{HHO18}, the set of all binomials corresponding to inner intervals of $J$ form a Gr\"{o}bner basis of $I_{J}$ under the order $<$.

The term order $<$ is the order $<^4$ defined in \cite{MRR22:primality_of_polyomino}.
 Since the minimal generating set of $I_{\MP}$ forms a reduced Gr\"{o}bner basis, $I_{\MP}$ is prime by~\cite[Corollary 3.3]{MRR22:primality_of_polyomino}.
\end{proof}

The polyominoes $\MP_1$ and $\MP_2$ shown in Figures~\ref{fig:examplethin} and~\ref{fig:examplethinpolyomino} satisfies the hypothesis of Theorem~\ref{thm:thingap2} respectively. Thus, for $i=1,2,$ $I_{\MP_i}$ is a Knutson prime ideal and the set of all binomials corresponding to inner intervals of $\MP_i$ forms a Gr\"{o}bner basis of $I_{\MP_i}$ under the order $<$.

\begin{figure}[t]{}
	\centering
	\begin{tikzpicture}[scale=.5]
		\draw[] (0,0) rectangle (11,1);
            \draw[] (6,1)--(6,2)--(7,2)--(7,1);
            \draw[] (9,0)--(9,-1)--(10,-1)--(10,0);
            \draw[] (10,0)--(10,2)--(11,2)--(11,1);
		\draw[] (1,-2) rectangle (2,10); 
		\draw[] (4,-1) rectangle (5, 6);
		\draw[] (-2,3) rectangle (12,4);
            \draw[] (-1,3)--(-1,2)--(-0,2)--(-0,3);
            \draw[] (-2,4)--(-2,5)--(-1,5)--(-1,4);
		\draw[] (8,0) rectangle (9,8);
		\draw[] (11,4) rectangle (12,2);
		\draw[] (1,7) rectangle (12,8);
            \draw[] (11,7)--(11,6)--(12,6)--(12,7);
            \draw[] (10,8)--(10,9)--(11,9)--(11,8);
            \draw[] (5,7)--(5,6)--(6,6)--(6,7);
		\draw[] (5,7) rectangle (6,9);
		\draw[] (1,-1) --(2,-1);
		\draw[] (1,2)--(2,2);
		\draw[] (1,5)--(2,5);
		\draw[] (1,6)--(2,6);
		\draw[] (1,7)--(2,7);
		\draw[] (1,9)--(2,9);
		\draw[] (8,2)--(9,2);
		\draw[] (8,5)--(9,5);
		\draw[] (8,6)--(9,6);
		\draw[] (3,0)--(3,1);
		\draw[] (6,0)--(6,1);
		\draw[] (7,0)--(7,1);
		\draw[] (3,3)--(3,4);
		\draw[] (0,3)--(0,4);
		\draw[] (-1,3)--(-1,4);
		\draw[] (6,3)--(6,4);
		\draw[] (7,3)--(7,4);
		\draw[] (10,3)--(10,4);
		\draw[] (3,7)--(3,8);
		\draw[] (4,7)--(4,8);
		\draw[] (7,7)--(7,8);
		\draw[] (10,7)--(10,8);
		\draw[] (11,7)--(11,8);
		\draw[] (4,5)--(5,5);
		\draw[] (4,2)--(5,2);
	\end{tikzpicture}
 \caption{}\label{fig:examplethinpolyomino}
\end{figure}

   \begin{figure}[h]
    \centering
		\begin{tikzpicture}[scale=.8]
  \draw[]
      (1,2)--(2,2)--(2,1)--(1,1)--(1,2)
        (2,3)--(3,3)--(3,2)--(2,2)--(2,3)
        (3,4)--(4,4)--(4,3)--(3,3)--(3,4)
        ;

        \filldraw[] (4,4) circle (1pt) node[anchor=west]  {$(i,j)$};
        \filldraw[] (3,3) circle (1pt) node[anchor=south east]  {$(i-1,j-1)$};
        \filldraw[] (2,2) circle (1pt) node[anchor=south east]  {$(i-2,j-2)$};
        \filldraw[] (1,1) circle (1pt) node[anchor=east]  {$(i-3,j-3)$};
	\end{tikzpicture}
    \caption{}\label{fig:corollary3}
 \end{figure}  

 \begin{remarkbox}
  Let $\MP$ be a thin collection of cells such that 
   \begin{asparaenum}
       \item
           it does not contain cells as in Figure~\ref{fig:corollary3};

       \item
            if $I$ is a maximal inner interval of $\MP$ with vertices $(i_1,j_1)$, $(i_1+1,j_1 )$, $(i_1,j_2 )$ and $(i_1+1, j_2 )$ (see Figure~\ref{fig:verticalinterval}), then the cell with upper-right vertex $(i_1,j_1)$ and the cell with lower-left vertex $(i_1+1,j_2)$ are not contained in $\MP$;

        \item
            if $I$ is a maximal inner interval of $\MP$ with vertices $(i_1,j_1)$, $(i_1,j_1+1)$, $(i_2,j_1 )$ and $(i_2, j_1+1)$ (see Figure~\ref{fig:horizontalinterval}), then the cell with upper-right vertex $(i_1,j_1)$ and the cell with lower-left vertex $(i_2,j_1+1)$ are not contained in $\MP$.
     \end{asparaenum}     
     Let $k\in \naturals$ be such that $(a,b)\notin V(\MP)$ for any $a, b\in \naturals$ with $a\geq k.$
     Let $\MQ$ be the collection of cells obtained by reflecting $\MP$ across the line $x=k$.
     Note that $\MQ$ satisfies the hypothesis of Theorem~\ref{thm:thingap2}.
     Also, one can see that there is an isomorphism $\varphi : S_\MP = K[x_v:v\in V(\MP)] \to S_\MQ=K[x_v:v\in V(\MQ)]$ defined as $\varphi(x_v = x_{(k-i,j)}) = x_{(k+i,j)}$. Then,
     $\varphi(I_\MP)=I_\MQ$. Therefore, we get that $I_\MP$ is a Knutson ideal and the set of all binomials corresponding to inner intervals of $\MP$ forms a Gr\"{o}bner basis of $I_{\MP}$. Moreover, when $\MP$ is a polyomino, then $I_\MP$ is a prime ideal.
 \end{remarkbox}

\begin{proposition}
    \label{prop:thinprime}
    Let $\MP$ be a thin polyomino such that if two maximal inner intervals of $\MP$ intersects, then there intersection is a cell.
    Then $I_\MP$ is a prime ideal.
\end{proposition}

\begin{proof}
Let $<$ be the reverse lexicographic order on $S$ induced by the following total order on variables of $S$: for $(i,j), (k,l)\in V(\MP)$, $x_{i,j}< x_{k,l}$, if $j<l$, or $l=j$ and $i<k$.
Under the hypothesis, by~\cite[Theorem\ 2.1]{JN24rookpolynomial}, the set of binomials corresponding to inner intervals of $\MP$ forms the reduced Gr\"{o}bner basis of $I_\MP$ with respect to $<$.
Also, note that the monomial order $<$ is the monomial order $<^6$ defined in~\cite{MRR22:primality_of_polyomino}.
Thus, by~\cite[Corollary\ 3.3]{MRR22:primality_of_polyomino}, $I_\MP$ is a prime ideal.
\end{proof}

\section{Some Knutson collection of cells}\label{sec:collectioncells}

In Section~\ref{sec:ladder}, we proved that the polyomino ideals associated with parallelogram polyominoes are Knutson. 
We now ask whether extracting a polyomino from a parallelogram polyomino is Knutson or not. 
In this section, we show that under certain hypothesis they are Knutson.

\begin{proposition}
        \label{prop:joiningpolyominoes}
       Let $\MP$ be a collection of cells such that $\MP=\MP_1\cup \MP_2$, where $\MP_1$ and $\MP_2$ are sub-collections of cells of $\MP$ with $\MP_1\cap \MP_2=\{v\}$, 
       a vertex of $\MP$. Suppose that  for each $i$, 
       $ I_{\MP_i}\in C_{f_i}$ for some $f_i\in S$ and $\gcd (x_v, \ini_<(f_i))=1$. 
        Then, there exists an $f\in S$ such that $I_{\MP} \in C_f$, that is, $I_{\MP}$ is a Knutson ideal.
\end{proposition}
\begin{proof}
    Take $f = f_1f_2x_v.$ By Lemma~\ref{subset}, $I_{\MP_i}\in C_f$ for $i=1,2.$ Now the result follows from the fact that $I_\MP = I_{\MP_1}+I_{\MP_2}$.
\end{proof}

Let $\MQ$ be a parallelogram polyomino.
Let $\MQ'$ be a sub-polyomino of $\MQ$ which is also a parallelogram polyomino.
Define $\MP = \MQ\setminus\MQ'$.
Note that $\MP$ may not be a polyomino.
Assume the notations of Discussion~\ref{setup:ladder} for the collection of cells $\MP$.
Set $a =\min\{i: V_i\cap V(\MQ')\neq \varnothing\}$ and $b=\max\{i: V_i\cap V(\MQ')\neq \varnothing\}$.
Define 
\[\MQ_1 = \MP_{b-1} = \cup_{i<b} \MC_i\]\
and 
\[\MQ_2 = \cup_{i>a} \MC_i. \]
Note that $\MQ_1$ and $\MQ_2$ may not be polyominoes, for example see Figure~\ref{fig:twosidedladderhole}.

 \begin{figure}[h]
    \centering
		\begin{tikzpicture}[scale=.8]
            \draw[] (2,1)--(9,1)--(9,2)--(9,3)--(10,3)--(10,4)--(10,5)--(11,5)--(11,6)--(6,6)--(6,5)--(4,5)--(4,3)--(3,3)--(3,2)--(2,2)--(2,1);
			\draw[] (3,2)--(9,2);
			\draw[] (3,3)--(6,3) (7,3)--(10,3);
			\draw[] (4,4)--(6,4) (8,4)--(10,4);
			\draw[] (5,5)--(11,5);
			\draw[] (6,6)--(8,6);
			\draw[] (2,1)--(2,2);
			\draw[] (3,1)--(3,2);
			\draw[] (4,1)--(4,5);
			\draw[] (5,1)--(5,5);
            \draw[] (6,1)--(6,2) (6,3)--(6,6);
			\draw[] (7,1)--(7,3) (7,5)--(7,6);
			\draw[] (8,1)--(8,6);
			\draw[] (9,1)--(9,6);
			\draw[] (10,5)--(10,6);
			\draw[] (11,5)--(11,6);

			\draw[color=blue] (4,3)--(6,1);

			\draw[color=blue] (7,6)--(10,3);
            
            \filldraw[] (4,3) circle (0pt) node[anchor=south east]  {$V_a$};
            \filldraw[] (7,6) circle (0pt) node[anchor=south east]  {$V_b$};

	\end{tikzpicture}
    \caption{}\label{fig:twosidedladderhole}
    \end{figure}
    
\begin{theorem}\label{thm:twosidedhole}
    Under the notations as above, if $I_{\MQ_1}+ I_{\MQ_2} = I_\MP$, then $I_\MP$ is Knutson.
\end{theorem}

\begin{proof}
We show that $I_\MP\in C_f$, where $f$ is as defined in Discussion~\ref{setup:ladder}.
It is enough to show that $I_{\MQ_1},I_{\MQ_2}\in C_f$. 
Let $h_1 = \cup_{i<b} f_i$ and  $h_2 = \cup_{i>a} f_i,$ where $f_i$ is as defined in Discussion~\ref{setup:ladder}.
By Lemma~\ref{subset}, it suffices to show that $I_{\MQ_i}\in C_{h_i}$ for all $i=1,2.$

We first assume that $\MP$ is non-simple.
We show that $I_{\MQ_1}\in C_{h_1}$. 
Notice that if $\MQ_1$ is a polyomino, then it is of the form as defined in Proposition~\ref{prop:ladderson interval}.
Otherwise, it is union of polyominoes $\MR, \MR_1^L,\ldots,\MR_r^L, \MR_1^R,\ldots,\MR_s^R$, where $\MR$ is a polyomino as defined in Proposition~\ref{prop:ladderson interval} and $\MR_i^L$'s and $\MR_j^R$'s are parallelogram polyominoes such that 
\begin{asparaenum}
\item 
$\MR\cap \MR_{1}^L$ (resp. $\MR\cap \MR_{1}^R$) is a vertex, say $v_1^L$ (resp. $v_1^R$); it is the minimal element of $\MR_{1}^L$ (resp. $\MR_{1}^R$) and it is the maximal element of $V(\MR_1)$ (resp. $V(\MR_2)$), where $\MR_1, \MR_2$ are as defined in Proposition~\ref{prop:ladderson interval} for $\MR$.
\item
$\MR_i^L\cap \MR_{i+1}^L$ is a vertex, say $v_{i+1}^L$; it is the maximal element of $V(\MR_i^L)$ and the minimal element of $V(\MR_{i+1}^L)$ for $1\leq i\ \leq r-1$.
\item 
$\MR_j^R\cap \MR_{j+1}^R$ is a vertex, say $v_{i+1}^R$; it is the maximal element of $V(\MR_j^R)$ and the minimal element of $V(\MR_{j+1}^R$) for $1\leq i\ \leq s-1$.
\end{asparaenum}
See Figure~\ref{fig:Q_1} for a schematic showing of the polyominoes $\MR, \MR_1^L,\ldots,\MR_r^L, \MR_1^R,\ldots,\MR_s^R$. 

If $\MQ_1$ is a polyomino, then we are done by Proposition~\ref{prop:ladderson interval}. So, we may assume that $\MQ_1$ is not a polyomino. 
Let $g, g_1^L,\ldots,g_r^L, g_1^R,\ldots,g_s^R$ be the polynomials as defined in Remark~\ref{rem:differentpolynomial} for $\MR, \MR_1^L,\ldots,\MR_r^L, \MR_1^R,\ldots,\MR_s^R$ respectively. Note that
\[
h_1 = x_v g x_{v_1^L}g_1^Lx_{v_2^L}g_2^L\cdots x_{v_r^L}g_r^L x_{v_1^R}g_1^Rx_{v_2^R}g_2^R\cdots x_{v_s^R}g_s^R.
\]
We have that $I_\MR\in C_g$ by Proposition~\ref{prop:ladderson interval}. 
Moreover, by Remark~\ref{rem:differentpolynomial}, $I_{\MR_i^L}\in C_{g_i^L}$ for all $1\leq i\leq r$ and $I_{\MR_j^R}\in C_{g_j^R}$ for all $1\leq i\leq s.$ We now use Proposition~\ref{prop:joiningpolyominoes} repeatedly to get that $I_{\MQ_1}\in C_{h_1}.$
Using a similar argument, one gets that $I_{\MQ_2}\in C_{h_2}.$

Consider the case when $\MP$ is simple. Note that if $\MQ_1$ is a polyomino, then we are done by Proposition~\ref{prop:ladderson interval} as before.
When $\MQ_1$ is not a polyomino, then it is a collection of cells as defined above or it is a collection of cells as defined above with $\MR,$ few of $\MR_i^L$'s and $\MR_j^R$'s are missing. In either case, by similar argument of non-simple case, one concludes that $I_{\MQ_1}\in C_{h_1}.$ Similarly, $I_{\MQ_2}\in C_{h_2}.$
\end{proof}

 \begin{figure}[t]
    \centering
		\begin{tikzpicture}[scale=.7]
            \draw[] (1,1)--(5,1)--(5,3)--(4,3)--(4,2)--(2,2)--(2,3)--(1,3)--(1,1)

            (2,3)--(2,5)--(3,5)--(3,3)--(2,3)
            (3,5)--(3,7)--(4,7)--(4,8)--(5,8)--(5,6)--(4,6)--(4,5)--(3,5)

           (6,10)--(8,10)--(8,9)--(6,9)--(6,10)


            (5,3)--(7,3)--(7,4)--(8,4)--(8,5)--(6,5)--(6,4)--(5,4)--(5,3)

            (9,6)--(10,6)--(10,8)--(9,8)--(9,6)
            

            ;
            \draw[dotted] 
                        (9,8)--(9,9)--(8,9)
                        (10,8)--(10,10)--(8,10)
                        (9,10)--(9,9)--(10,9)
            ;
            \draw[color=blue] (8,10)--(10,8);
            
            \filldraw[] (3,1.5) circle (0pt) node[anchor=center]  {$\MR$};
            \filldraw[] (6.5,4) circle (0pt) node[anchor=center]  {$\MR_1^R$};
            \filldraw[] (9.5,7) circle (0pt) node[anchor=center]  {$\MR_s^R$};
            \filldraw[] (8.5,5.5) circle (0pt) node[anchor=center]  {$\iddots$};
            
            \filldraw[] (2.5,4) circle (0pt) node[anchor=center]  {$\MR_1^L$};
            \filldraw[] (4,6.5) circle (0pt) node[anchor=center]  {$\MR_2^L$};
            \filldraw[] (7,9.5) circle (0pt) node[anchor=center]  {$\MR_r^L$};
            \filldraw[] (5.5,8.5) circle (0pt) node[anchor=center]  {$\iddots$};

            \filldraw[] (8,10) circle (0pt) node[anchor=south east]  {$V_b$};

            \filldraw[] (1,1) circle (1pt) node[anchor=east]  {$v$};

            \filldraw[] (5,3) circle (1pt) node[anchor=north west]  {$v_1^R$};
            \filldraw[] (8,5) circle (1pt) node[anchor=north west]  {$v_2^R$};
            \filldraw[] (9,6) circle (1pt) node[anchor=north west]  {$v_s^R$};
            
            \filldraw[] (2,3) circle (1pt) node[anchor=south east]  {$v_1^L$};
            \filldraw[] (3,5) circle (1pt) node[anchor=south east]  {$v_2^L$};
            \filldraw[] (5,8) circle (1pt) node[anchor=south east]  {$v_3^L$};
            \filldraw[] (6,9) circle (1pt) node[anchor=south east]  {$v_r^L$};

	\end{tikzpicture}
    \caption{}\label{fig:Q_1}
    \end{figure}

We now give a characterization when $I_{\MQ_1}+ I_{\MQ_2} = I_\MP$.

\begin{proposition}
Under the notations as above, $I_{\MQ_1}+ I_{\MQ_2} = I_\MP$ if and only if there is no inner interval of $\MP$ that contain cells from both $\MC_a$ and $\MC_b$. 
\end{proposition}   

\begin{proof}
Since $\MQ_1$ and $\MQ_2$ are sub-collections of cells of $\MP,$ $I_{\MQ_1}+ I_{\MQ_2} \subseteq I_\MP.$
Also, note that by definition of $\MQ_1$ and $\MQ_2$, no cell from $\MC_b$ is included in $\MQ_1$, and no cell from $\MC_a$ is included in $\MQ_2$.

    $(\implies)$ On the contrary, assume that $J$ is an inner interval of $\MP$ that contains cells from both $\MC_a$ and $\MC_b$.
    This means that $J$ cannot be an inner interval of either $\MQ_1$ or $\MQ_2$. Therefore, the generator of $I_\MP$ corresponding to $J$ does not belong to $I_{\MQ_1}$ and $I_{\MQ_2}$. Then, one can deduce that, $I_\MP \not\subset I_{\MQ_1}+ I_{\MQ_2}$ which is a contradiction.
    
    $(\impliedby)$ It suffices to show that every inner interval of $\MP$ is also an inner interval of either $\MQ_1$ or $\MQ_2.$ 
    On the contrary, assume that $I$ is an inner interval of $\MP$ that is not an inner interval of either $\MQ_1$ or $\MQ_2.$ This implies that $I$ contains cells from $\MC_i$ for some $i\leq a$ and from $\MC_j$ for some $j\geq b.$ Since $I$ is an inner interval and, by definition of $V_i$'s, it must also contain cells from both $\MC_a$ and $\MC_b$. This leads to a contradiction.
\end{proof}

\begin{proposition}
    Let $\MP, \MQ_1$ and $\MQ_2$ be as above such that $I_{\MQ_1}+ I_{\MQ_2} = I_\MP$. Then, the set of all quadratic binomials corresponding to inner intervals of $\MP$ form a Gr\"{o}bner basis of $I_{\MP}$ under the order $<$ as defined in Discussion~\ref{setup:ladder}. 
\end{proposition}
\begin{proof}
    Since $I_{\MQ_1}, I_{\MQ_1}\in C_f$ by Theorem~\ref{thm:twosidedhole}, and $I_{\MQ_1}+ I_{\MQ_2} = I_\MP$, it suffices to prove the result for $\MQ_1$ and $\MQ_2$. We first proceed for $\MQ_1.$
 
    When $\MQ_1$ is a polyomino, it is as defined in Proposition~\ref{prop:ladderson interval} such that $\MR_1$ and $\MR_2$ are two polyominoes based on a polyomino $\MI$, where $\MI, \MR_1$ and $\MR_2$ are as defined in Proposition~\ref{prop:ladderson interval}. Since $\MP$ is obtained by removing a parallelogram polyomino from a parallelogram polyomino, the polyominoes $\MI, \MR_1$ and $\MR_2$ are also parallelogram polyominoes. So we are done by Proposition~\ref{prop:groebner}.

    On the other hand, when $\MQ_1$ is not a polyomino, it is a collection of cells as defined in the proof of Theorem~\ref{thm:twosidedhole}. Assume the notations of proof of Theorem~\ref{thm:twosidedhole}. Note that $\MR$ is a polyomino as described in the previous case. Moreover, Since $\MR_i^L$'s and $\MR_j^R$'s are parallelogram polyominoes, the $K$-algebras $K[\MR_i^L]$'s and $K[\MR_j^R]$'s are the Hibi ring associated to the vertex sets  $V(\MR_i^L)$'s and $V(\MR_j^R)$ which are distributive lattices (see~\cite{KV21}, \cite[Remark\ 2.4]{QRR22}). By~\cite[Theorem\ 6.17]{HHO18}, the set of all binomials corresponding to inner intervals of $\MR_i^L$'s and $\MR_j^R$'s form a Gr\"{o}bner basis of $I_{\MR_i^L}$'s and $I_{\MR_j^R}$'s under the order $<$. Since $I_{\MQ_1} = I_\MR + \sum_{i=1}^{r} I_{\MR_i^L} +\sum_{j=1}^{s} I_{\MR_j^R},$ the result follows for $\MQ_1$.

    Using a similar argument, one get that the set of all quadratic binomials corresponding to inner intervals of $\MQ_2$ form a Gr\"{o}bner basis of $I_{\MQ_2}$ under the order $<$
\end{proof}


\newcommand{\etalchar}[1]{$^{#1}$}

\end{document}